\documentclass[12pt]{amsart}
\usepackage{amsfonts,amssymb,amscd,amsmath}
\usepackage{latexsym}
\usepackage{tocvsec2}

\usepackage[
hyperindex=true,pagebackref=true,bookmarks=true,
colorlinks=true,linkcolor=blue,citecolor=red]
{hyperref}
\def\~{{\rm --}}

\begin{document}
\renewcommand{\theequation}{\thesection.\arabic{equation}}
\settocdepth{subsubsection}
\textheight 7.5truein

\newcommand{\Z}{{\mathbb Z}}
\renewcommand{\S}{{\mathbb S}}
\newcommand{\Q}{{\mathbb Q}}
\newcommand{\N}{{\mathbb N}}
\newcommand{\C}{{\mathbb C}}
\newcommand{\R}{{\mathbb R}}
\newcommand{\T}{{\mathbb T}}
\newcommand{\G}{{\mathbb G}}
\newcommand{\E}{{\mathbb E}}
\newcommand{\K}{{\mathbb K}}
\newcommand{\msl}{\mathfrak{sl}}
\newcommand{\mgl}{\mathfrak{gl}}
\newcommand{\mso}{\mathfrak{so}}
\newcommand{\msp}{\mathfrak{sp}}
\newcommand{\veps}{\varepsilon}
\newcommand{\fg}{{\mathfrak g}}
\newcommand{\fb}{{\mathfrak b}}
\newcommand{\fn}{{\mathfrak n}}
\newcommand{\fh}{{\mathfrak h}}

\newtheorem{theorem}{Theorem}[section]
\newtheorem{conjecture}[theorem]{Conjecture}
\newtheorem{maintheorem}[theorem]{Main Theorem}
\newtheorem{proposition}[theorem]{Proposition}
\newtheorem{definition}[theorem]{Definition}
\newtheorem{lemma}[theorem]{Lemma}
\newtheorem{corollary}[theorem]{Corollary}
\newtheorem{notation}[theorem]{Notation}
\newtheorem{remark}[theorem]{Remark}
\newtheorem{example}[theorem]{Example}

\def\equal{\stackrel{\,\mathbf{def}}{= \kern-3pt =}}

\def\la{\lambda}
\def\La{\Lambda}
\def\om{\omega}
\def\Om{\Omega}
\def\Th{\Theta}
\def\th{\theta}
\def\al{\alpha}
\def\tal{\tilde{\alpha}}
\def\be{\beta}
\def\tbe{\tilde{\beta}}
\def\ga{\gamma}
\def\tga{\tilde{\gamma}}
\def\ep{\epsilon}
\def\up{\upsilon}
\def\Up{\Upsilon}
\def\de{\delta}
\def\De{\Delta}
\def\ka{\kappa}
\def\kapp{\hbox{\bf \ae}}
\def\si{\sigma}
\def\Si{\Sigma}
\def\Ga{\Gamma}
\def\ze{\zeta}
\def\io{\iota}
\def\bep{\backepsilon}
\def\pa{\partial}
\def\vph{\varphi}
\def\vep{\varepsilon}
\def\vpi{{\varpi}}
\def\vth{{\vartheta}}
\def\vsi{{\varsigma}}
\def\vrh{{\varrho}}

\def\tR{\tilde{R}}
\def\hw{\hat{w}}
\def\hW{\hat{W}}

\def\f{\mathcal{F}}
\def\o{\mathcal{O}}
\def\t{\mathcal{T}}
\def\r{\mathcal{R}}
\def\l{\mathcal{L}}
\def\m{\mathcal{M}}
\def\k{\mathcal{K}}
\def\n{\mathcal{N}}
\def\d{\mathcal{D}}
\def\p{\mathcal{P}}
\def\a{\mathcal{A}}
\def\h{\mathcal{H}}
\def\c{\mathcal{C}}
\def\y{\mathcal{Y}}
\def\e{\mathcal{E}}
\def\v{\mathcal{V}}
\def\z{\mathcal{Z}}
\def\x{\mathcal{X}}
\def\s{\mathcal{S}}
\def\g{\mathcal{G}}
\def\u{\mathcal{U}}
\def\w{\mathcal{W}}
\def\i{\mathcal{I}}
\def\j{\mathcal{J}}
\def\b{\mathcal{B}}

\def\lng{\hbox{\rm{\tiny lng}}}
\def\sht{\hbox{\rm{\tiny sht}}}
\def\sing{\hbox{\rm{\tiny sing}}}
\def\lan{\langle}
\def\ran{\rangle}
\font\smm=msbm10 at 12pt 
\def\symbol#1{\hbox{\smm #1}}
\def\lsmash{{\symbol n}}
\def\rsmash{{\symbol o}}
\def\#{\sharp}
\newcommand{\sq}{\phantom{1}\hfill$\qed$}

\renewcommand{\tilde}{\widetilde}
\renewcommand{\hat}{\widehat}
\newcommand{\comment}[1]{}

\title
[Extremal part of the PBW-filtration and {\em E\,}-polynomials]
{Extremal part of the PBW-filtration and {\em E\,}-polynomials}
\author
[Ivan Cherednik and Evgeny Feigin]
{Ivan Cherednik and Evgeny Feigin}

\begin{abstract}
Given a reduced irreducible root system, the corresponding nil-DAHA 
is used to calculate the extremal coefficients of nonsymmetric 
Macdonald polynomials, also called $E\,$-polynomails, in the limit 
$t\to \infty$ and for antidominant weights, which is an important 
ingredient of the new theory of nonsymmetric $q$-Whittaker function. 
These coefficients are pure $q$-powers and their degrees are expected 
to coincide in the untwisted setting with the extremal degrees of 
the so-called PBW-filtration in the corresponding finite-dimensional 
irreducible representations of the simple Lie algebras for any root 
systems. This is a particular case of a general conjecture in terms 
of the level-one Demazure modules. We prove this coincidence for all 
Lie algebras of classical type and for $G_2$, and also establish the 
relations of our extremal degrees to minimal $q$-degrees of the 
extremal terms of the Kostant $q$-partition function; they coincide 
with the latter only for some root systems. 
\end{abstract}

\thanks{ \today}
\address[I. Cherednik]{Department of Mathematics, UNC
Chapel Hill, North Carolina 27599, USA\newline
chered@email.unc.edu}
\address[E. Feigin]{Department of Mathematics, 
National Research University Higher School of Economics,
Vavilova str. 7, 117312, Moscow, Russia,
{\it and }
Tamm Theory Division, Lebedev Physics Institute\newline
evgfeig@gmail.com}

\maketitle
\smallskip

{\small
\centerline{{\bf MSC} (2010): 17B45, 22C08, 33D52, 17B67}
}

\renewcommand{\baselinestretch}{1.2}
{\textmd
\tableofcontents
}
\renewcommand{\baselinestretch}{1.0}
\vfill\eject

\setcounter{section}{0}
\setcounter{equation}{0}
\section{\sc Introduction}
The nil-DAHA, more specifically the theory of
the so-called $E$-dag polynomials, is employed in 
this paper to obtain new surprising formulas for 
the extremal degrees of the PBW-filtration 
in finite-dimensional irreducible representations 
of simple Lie algebras of classical type and $G_2$ 
in the untwisted setting. This is expected to hold
for any untwisted (reduced) root systems with certain 
expectations in the twisted setting. This correspondence
is a particular case of a general conjecture in terms 
of the level-one Demazure modules. 

The key step here is in establishing the additivity of the 
formulas for the extremal dag-degrees and for the
extremal PBW-degrees in the (anti)dominant sector
(dominant in PBW-theory and antidominant for the
$E$-dag polynomials); then the fundamental weights are 
sufficient to consider. The extremal dag-degrees for the 
latter are provided for $A\!BC\!D\!FG$, to be
systematically considered in further works.   

We also discuss the relation of our formulas to 
the minimal $q$-degrees of the extremal part of 
the Kostant $q$-partition function. This function 
is connected with both theories, nil-DAHA and PBW,
but its minimal $q$-degrees coincide with ours for all
(anti)dominant weights only for types $A$, twisted $B$, 
untwisted $C$ and twisted $G_2$. The extremal part of
the Kostant $q$-partition function is generally not
additive.

Our result is a special case of a general conjecture
connecting full $E$-dag polynomials with
the Demazure level-one modules supplied with
the sum of the Kac-Moody-degree and the PBW-degree.
Upon the restriction to the $W$-extremal vectors, the 
Kac-Moody filtration vanishes and we obtain a
surprising  application of DAHA, double affine Hecke algebras,
to the classical theory of finite-dimensional representations 
of simple Lie algebras, which seems the first 
``nonaffine" DAHA application of this scale. However our proof 
is of technical nature and does not clarify the real reasons 
of this correspondence (and those behind the general 
affine conjecture).
\smallskip

\subsection{\bf \texorpdfstring{$E$}{E}-polynomials and 
\texorpdfstring{$E$}{E}-dag polynomials}
The dag-polynomials are the limits of
the nonsymmetric Macdonald polynomials as $t=\infty$;
they are dual to the generalized nonsymmetric
$q$-Hermite polynomials (corresponding to $t\to 0$), called 
for short $E$-bar polynomials in this paper. The duality is with 
respect to the inner product in terms
of the standard  multiplicative theta-function associated with 
a given root system. 
See \cite{Op, Ma, Ch1, Ch2} for general 
theory of nonsymmetric Macdonald polynomials, also
called $E$-polynomials. 

The formulas for the extremal $q$-degrees of $E$-dag polynomials 
are presented in this paper for classical root systems and 
$G_2$, as well as computer-generated formulas for  $F_4$ and 
for $E_6$ (provided only partially; they were 
calculated for $E_7$ too). For the classical root systems, 
the extremal degrees for the fundamental weights can be calculated
by a relatively straightforward induction.  
The $q$-positivity of complete $E$-dag polynomials was conjectured in 
\cite{CO1}, which is a theorem for their extremal parts and for 
antidominant weights (see Theorem \ref{EXTDAG}\, below;
it was announced in \cite{CO1}, Corollary 2.6).


The formulas for extremal $q$-degrees 
of $E$-dag polynomials for {\em all} weights are a significant 
ingredient of the new theory of nonsymmetric $q$-Whittaker 
function \cite{CO1}. This link is expected to be important to
understand their meaning, but we present (and partially justify)
the {\em PBW-$E^\dag$ correspondence\,} in an entirely algebraic way 
in this paper.
\medskip

\subsection{\bf Hall-Littlewood and \texorpdfstring{$E$}{E}-bar 
polynomials}
An important development of the classical theory
of finite-dimensional representations of semisimple 
finite-dimensional Lie algebras was the introduction
of the so-called BK-filtration (see \cite{Ko,Br,JLZ}) and 
establishing its relation with Lusztig's $q$-analogs 
of weight multiplicities, defined in \cite{Lu} via the
affine Kazhdan-Lusztig polynomials upon their restriction
to the lattice of radical weights. This theory is
directly connected with the Hall-Littlewood polynomials, 
equivalently, Macdonald's $p$-adic spherical functions. 
We use this theory as a natural pattern, 
but the PBW-filtration and the BK-filtration are really 
different, as well as the corresponding polynomials. 

The Hall-Littlewood polynomials are the limit $q\to 0$ 
of the Macdonald symmetric polynomials. Such limit  
(among other simplifications) results in explicit formulas 
for these polynomials, instead of obtaining them as eigenfunctions 
of certain $q$-difference operators in the general theory.
The $E$-bar polynomials ($t\to 0$) have important
applications too. They coincide with the level-one Demazure 
characters in the twisted setting \cite{San,Ion1} for all weights, 
not only (anti)dominant. They are also related to the characters
of the local and global Weyl modules (see, e.g. \cite{FeL},
\cite{FoL} and \cite{CL}).
Let us mention here well-known and fruitful relations of the 
Schubert polynomials and quantum Schubert polynomials
to nil-Hecke algebras and similar objects; see \cite{FGP} and 
references therein. This can be connected with our usage 
of nil-DAHA. 

Also, the $E$-bar polynomials
and the {\em global} (symmetric) $q$-Whittaker function,
which is a quadratic generating function of the $E$-bar polynomials
for antidominant weights, are directly related to the
Gromov-Witten invariants of flag varieties 
and affine flag varieties; see \cite{GL} and \cite{BF}. 
The main link is via the Harish-Chandra-type asymptotic
expansions of global $q$-Whittaker functions, but there
are other important aspects of this relation. 
The global $q$-Whittaker and $q$-hypergeometric functions 
are actually of algebraic nature (in contrast to those 
without $q\,$ in the
classical harmonic analysis on symmetric spaces).
\medskip

\subsection{\bf PBW-filtration and \texorpdfstring{$E$}{E}-dag 
polynomials}
The main result of this paper is in establishing
the correspondence between the extremal $E$-dag 
polynomials and the PBW-filtration. This is the 
extremal part of Conjecture 2.7 from \cite{CO1} (with
participation of E.~F.) on the coincidence of the 
$E$-dag polynomials and the corresponding 
characters of Demazure level-one modules for
the sum of the Kac-Moody grading and that 
defined via the PBW-filtration. 
It was/is stated only for (anti)dominant weights, which is
significant, and for the $A\!D\!E$ root systems; the 
twisted case is in progress. 

Calculating the PBW-filtration 
is generally a difficult problem; no systematic methods beyond 
(nonaffine) types $A$, untwisted $C$ and untwisted $G_2$
are known at the moment 
\cite{F1,FFL1,FFL2,FFL3,Gor}. This becomes especially involved in the 
Kac-Moody case (the setting of Conjecture 2.7 from \cite{CO1}
and Conjecture \ref{CONJDEM} below).
The extremal coefficients of the $E$-dag polynomials for antidominant 
weights can be calculated for all root systems
(using computers for $F_4,E_{6,7,8}$), so the problem
with their identification is currently due to the lack of methods 
on the PBW-side of this surprising correspondence. 
\smallskip

The twisted setting is not a problem for the $E$-dag polynomials
(it is actually a preferred setup). However by now there is no
twisted PBW-theory. 
Actually, the {\em PBW-$E^\dag$ correspondence\,}
is quite a challenge even when it is justified (untwisted classical
root systems and $G_2$ for the extremal terms).   
We hope that this correspondence is of geometric nature 
(to be discovered) similar to the theory of $E$-dag polynomials 
the Hall-Littlewood polynomials, though the latter two families
seem significantly simpler than the $E$-dag one. 

\medskip

\setcounter{equation}{0}
\section{\sc Affine root systems}
\subsection{\bf Basic notations}
Let $R=\{\al\} \subset \R^n$ be a root system of type
$A,B,...,F,G$
with respect to a Euclidean form $(z,z')$ on $\R^n
\ni z,z'$,
$W$ the {\em Weyl group}
generated by the reflections $s_\al$,
$R_{+}$ the set of positive  roots ($R_-=-R_+$)
corresponding to fixed simple
roots $\al_1,...,\al_n,$
$\Ga$ the Dynkin diagram
with $\{\al_i , 1 \le i \le n\}$ as the vertices.
Accordingly,
$$R^\vee=\{\al^\vee =2\al/(\al,\al)\}.$$

The root lattice and the weight lattice are:
\begin{align}
&Q=\oplus^n_{i=1}\Z \al_i \subset P=\oplus^n_{i=1}\Z \om_i,
\notag
\end{align}
where $\{\om_i\}$ are fundamental weights:
$(\om_i,\al_j^\vee)=\de_{ij}$ for the
simple coroots $\al_i^\vee.$
Replacing $\Z$ by $\Z_{\pm}=\{m\in\Z, \pm m\ge 0\}$ we obtain
$Q_\pm, P_\pm.$
Here and further see \cite{Bo}.

The form will be normalized
by the condition  $(\al,\al)=2$ for 
{\em short} roots in this paper.
The normalization leads to the inclusions
$Q\subset Q^\vee,  P\subset P^\vee,$ where $P^\vee$ is
generated by the fundamental coweights $\{\om_i^\vee\}$
dual to $\{\al_i\}$. We set 
\begin{align}\label{etadef}
\nu_\al\ =\ (\al,\al)/2,\ \, \nu_i=\nu_{\al_i} \hbox{\  for\  }
1\le i \le n.
\end{align} 

We will constantly use  
\begin{align}\label{partialrho}
&\rho\,\equal\, \frac{1}{2}\sum_{\al\in R_+}\,\al=
\sum_{i=1}^n \om_i,\ \ 
\check{\rho}\,\equal\, \frac{1}{2}\sum_{\al\in R_+}\,\al^\vee=
\sum_{i=1}^n \om_i^\vee.
\end{align}

Two maximal roots will be considered in this paper,
the standard maximal positive root $\th_{\lng}=\th$ and
the maximal short root $\th_{\sht}=\vth\in R_+$. 
The latter is the {\em maximal positive
coroot} because of the choice of normalization.

\smallskip

\subsection{\bf Affine root systems}
The {\em affine untwisted root system} is
$$
\tilde{R}=\{\tilde{\al}=[\al,j]\,:\, \al\in R,\, j\in \Z\}\,
\subset\, \R^n\times \R \subset \R^{n+1}.
$$
We identify $z\in \R^n$ with $ [z,0]$, so $R\subset \tilde{R}$.
Accordingly $\tilde{R}_\pm=\{\tilde{\al}\,:\, j>0 \hbox{ or }
j=0, \al>0\}$. We add $\al_0=[-\th,1]$ to the set of simple
roots and denote the completed Dynkin diagram 
by $\tilde{\Ga}$. 

The other affine extension of $R$ is the
{\em twisted affine root system} $\tilde{R}^\nu$ formed by
vectors $\ \tilde{\al}=[\al,\nu_\al j]\ $
for $\al \in R, j \in \Z $. The corresponding $\al_0$ 
is $\al^\nu_0 \equal [-\vth,1]$ for $\vth=\th_{\sht}\in R_+$.
Setting $\tal^\vee=\tal/\nu_\al$ affine roots,
$\tR^\nu=(\widetilde{R^\vee})^\vee$.
We will frequently omit the super-index $\nu$ in the twisted
case using the same notation $\tilde{R}$ and $\al_0$ as in the 
untwisted case, unless misunderstanding is not impossible. 
Also, the notation $\tilde{R}^\backepsilon$ will be frequently
used for the sake of uniformity in the formulas stated for 
both settings.

The {\em twisted completed Dynkin diagram}
$\tilde{\Ga}^\nu$ is obtained from  $\Ga$ 
by adding short $\al_0^\nu$ ($-\vth$, to be more
exact). One can obtain it from the
completed Dynkin diagram from \cite{Bo} for
the dual system
$R^\vee$ by reversing all arrows. 

The set of indices of the images of $\al_0$ by all
the automorphisms of $\tilde{\Ga}$ will be denoted by $O$
($O=\{0\} \hbox{\ for\ } E_8,F_4,G_2$); \, $O'\equal
O\setminus \{0\}$. We will use the
same notation $O$ for the orbit of $\al_0^\nu$
in $\tilde{\Ga}^\nu$. In the twisted setting
the elements $\om_r$ for $r\in O'$ are minuscule
weights: $(\om_r,\al^\vee)\le 1$ for all $\al \in R_+$.
We set here and below $\om_0=0$. 
\smallskip

\subsection{\bf Affine Weyl groups}
In the twisted or untwisted case, they are 
generated by all $s_{\tilde{\al}}$
for $\tilde{\al}\in\tilde{R}_+$;
we write $\tilde{W} = \lan s_{\tilde{\al}},\,\tilde{\al}\in 
\tilde{R}_+\ran$. One can take
the simple reflections $s_i=s_{\al_i}\ (0 \le i \le n)$
as its
generators and introduce the corresponding notion of the
length (see below). Note that the only difference between 
the twisted and untwisted cases is in the definition of $\al_0$. 
We will add the super-index $\nu$ to emphasize (when necessary) 
that the twisted case is considered;
for instance, $\tilde{W}^\nu=\lan s_{\tal},\tal\in 
\tR_+^\nu\ran$ denotes the twisted affine Weyl group. 

Given $\tilde{\al}=[\al,j]$ and $a \in P$, 
\begin{align}
&s_{\tilde{\al}}(\tilde{z})\ =\ \tilde{z}-(z,\al^\vee)\tilde{\al},\
\ a(\tilde{z})\ =\ [z,\ze-(z,a)]
\label{ondon}
\end{align}
for $\tilde{z}=[z,\ze] \in \R^{n+1}$. These formulas do not
depend on the setting, twisted or untwisted. We extend
the form $(\cdot,\cdot)$ to $\R^{n+1}$ by 
$([z,\ze],[z',\ze'])=(z,z')$ and set $\tal^\vee=\tal/\nu_\al$.
Thus one can use the inner product 
$(\tilde{z},\tal^\vee)$ in (\ref{ondon}) instead of $(z,\al^\vee)$.
 
The group $\tilde{W}$ is the semidirect product 
$W\lsmash Q^\vee$ of
its subgroups $W=\lan s_\al, \al \in R_+\ran$ and $Q^\vee$
in the untwisted setting and 
$\tilde{W}^\nu=W\lsmash Q$ in the twisted setting, where
the elements of $Q$ and $Q^\vee$ act in $R^{n+1}$ via
the second formula in (\ref{ondon}).
For instance for $\al\in R$,
\begin{align}
& Q^\vee\!\ni\! \al^\vee\!=\!s_{\al}s_{[\al,\,1]}=
s_{[-\al,\,1]}s_{\al},\ \,
Q\!\ni\! \al\!=\!s_{\al}s_{[\al,\,\nu_{\al}]}=
s_{[-\al,\,\nu_\al]}s_{\al}.
\label{ondtwo}
\end{align}

By $l(\tilde{w})$, we mean  the length of the minimal (reduced)
decomposition of $\tilde{w}$ in terms of simple reflections.
It can be also defined as the
cardinality $\#\{\la(\tilde{w})\}$
of the {\em \,$\la$\~set\,} of $\tilde{w}$\,:
\begin{align}\label{lasetdef}
&\la(\tilde{w})\equal \tilde{R}_+\cap 
\tilde{w}^{-1}(\tilde{R}_-)=\{\tilde{\al}\in R_+,\
\tilde{w}(\tilde{\al})\in \tilde{R}_-\},\
\tilde{w}\in \tilde{W}.
\end{align}
One has
\begin{align}
& \la(\tilde{w}\tilde{u}) = \la(\tilde{u}) \cup
\tilde{u}^{-1}(\la(\tilde{w})) \hbox{\ \,provided\ \,} 
l(\tilde{w}\tilde{u})=l(\tilde{w})+l(\tilde{u}).
\label{ltutw}
\end{align}
In the twisted case, $\tilde{R}^\nu$ must be
used instead of $\tilde{R}$
in the definition of the $\la$-sets. We will use
that $l(b)=2(\check{\rho},b)$ for $b\in P$ and 
$l(b)=2(\rho,b)$ for $b\in P^\vee$. Using the uniform
notation, $l(b)=2(\check{\rho}^\bep,b)$, where
$\bep\,=\varnothing$ in the untwisted case and
$\bep\,= \nu$ in the twisted case; $\check{\rho}^\nu=\rho$.

\smallskip

\subsection{\bf Extended Weyl groups} We define 
$\hat{W}$ and $\hat{W}^\nu$, correspondingly, 
as $W\lsmash P^\vee$ and $W\lsmash P$ acting in 
$\R^{n+1}$ via (\ref{ondon}):
\begin{align}
&(wa)([z,\ze])\ =\ [w(z),\ze-(z,a)] \for w\in W,\ a\in \,P^\vee,P.
\label{ondthr}
\end{align}

Given $a\in\, P_+^\vee,P_+$, let $w^a_0$ be the longest element
in the subgroup $W^{a}\subset W$ of the elements
preserving $a$. This subgroup is generated by simple
reflections. We set
\begin{align}
&u_{a} = w_0w^a_0  \in  W,\ \pi_{a} =
a( u_{a})^{-1} \in \,\hat{W},\hat{W}^\nu, \  
u_i= u_{\om_i},\ \pi_i=\pi_{\om_i},
\label{xwo}
\end{align}
where $w_0$ is the longest element in $W,$
$1\le i\le n.$ More generally, $u_a$ is the 
greatest element from $W$ such that the decomposition
$a=\pi_a u_a$ is reduced, i.e. $l(a)=l(\pi_a)+l(u_a)$;
it can be defined in this way for any $\hat{w}\in \hat{W}$
instead of $a\in P_+$.

Recall that $O$ is the orbit of $\al_0$ or $\al_0^\nu$
in $\tilde{\Ga}$ or $\tilde{\Ga}^\nu$
under the action of the group of its automorphisms.
Also, $O' \equal O\setminus \{0\}$ and $\pi_0=\hbox{id}$;
in the twisted setting $O'$ is the set of indices of 
minuscule weights (coweights in the untwisted case).  
The elements $\pi_r=\pi_{\om_r}$ for the indices
$r\in O$ leave $\tilde{\Ga}$ invariant. 
They form a group denoted by $\Pi$,
which is isomorphic to $P/Q$ by the natural
projection $\{\om_r \mapsto \pi_r\}$ in the twisted
case and by $\{\om_r^\vee \mapsto \pi_r\}$ in the
untwisted case. 

We set $l(\pi_r)=0$, extending (\ref{lasetdef})
to $\hat{W}$. Switching to $u_r \ (r\in O')$,
these elements preserve the set  
$\{-\th_{\sht,\lng}\}\cup\{\al_i, i>0\}$;
recall that $\vth=\th_{\sht}$ and $\th=\th_{\lng}$ 
are taken in the twisted/untwisted cases.

The relations $\pi_r(\al_0)= \al_r$
distinguish the indices $r \in O$. Moreover, one has
\begin{align}
& \hat{W}^\backepsilon  = \Pi \lsmash \tilde{W}^\backepsilon, 
\hbox{\ where\ }
\pi_rs_i\pi_r^{-1}  =  s_j \hbox{\ \ iff\ \ }
\pi_r(\al_i)=\al_j,\ i,j\ge 0.
\end{align}
Here backepsilon $\bep\,$ is $\nu$ in the twisted case 
or $\varnothing$ in the untwisted case.
\smallskip

\subsection{\bf Classical Weyl groups}
For $A_n$, the nonaffine Weyl group is
$W=\S_{n+1}$. The standard one-line
notation $w=(w_1,\cdots w_{n+1})$ will be used for permutations;
$s_{ij}=(ij)$ is the transposition of $i$ and $j$,\,
$s_i=s_{i\,i+1}$ for $1\le i\le n$. We will also switch
to $w(k)=w_k$ when it is convenient.

For the other classical root systems, using the
permutations with signs is standard. We represent
$W=\{w=(w_i)\,:\, 1\le i\le n\}$, where $(|w_i|)$ is
a permutation of the set of $\{1,2,\ldots,n\}$. The 
signs of $w_i$ can be arbitrary for $B_n,C_n$ and the total 
of negative $w_i$ must be even for $D_n(n>3)$. They are 
composed naturally; namely, we interpret such $w$ as 
transformations
$$\{\ep i\mapsto \ep w_i\,:\, 1\le i\le n\}
\hbox{\, of the set\, }\{\pm i:\, 1\le i\le n\}, \hbox{\, where }
\ep=\pm 1.
$$
For instance, $w_0=(n+1,\ldots,1)$ for $A_n$,
$w_0=(-1,\ldots,2-n,-n,1-n)$ for $D_n$ with odd $n$,
and $w_0=(-1,-2,\ldots,-n)$ otherwise.

The simple reflections are 
$s_i$ for $i<n$ for $A_n$,
$s_n=(1,2,\ldots,n-1,-n)$ for $B_n,C_n$ and
$s_n=(1,2,\ldots,n-2,-n,1-n)$ for $D_n(n>3)$.

We will also use the formulas for the fundamental
weights and $\th_{\lng,\sht}$ following the notation
from \cite{Bo}. Recall that $\om_1,\om_n,\om_{n-1}$ are
minuscule for $D_n(n\ge 4)$ and all fundamental roots are
such for $A_n$.  For $B_n,C_n$, the set $O'$ is
correspondingly  $\{1\}$, $\{n\}$ in the untwisted case
and $\{n\}$, $\{1\}$ for the twisted setting.

For $1\le i\le n$, let us list the stabilizers $W^i=W^{\om_i}$:
\begin{align}\label{stabwi}
A_n&,B_n,C_n\,:\, \{(w_j\,\mid\, 0<w_j\le i \hbox{\, for\, } j\le i\},
\\
D_n&\,:\, \{(w_j\,\mid\, 0<w_j\le i \hbox{ for } j\le i\neq n-1\},
\notag\\
 &\,:\,\{(\ze_j|w_j| \hbox{\, for\, } j\le i+1=n\},
\hbox{\,\ where\,\ }
 \ze_j=1\notag\\ 
&\hbox{\,\ \ \ \ \ except for \,} \ze_j\!=\!-1 \hbox{\, when } j=n 
\hbox{\, or } |w_j|=n.\notag
\end{align}



\setcounter{equation}{0}
\section{\sc Extremal dag-polynomials}
Let $\Z_q[X_b, b\in P]$ be the algebra of
Laurent polynomials in terms of $X_b$ satisfying
$X_{b+c}=X_b X_c$ with the coefficients in
$\Z_q\equal\Z[q,q^{-1}]$. The construction below
is for formal $q$, but one can substitute any $q\neq 0$
(including roots of unity).
In contrast to the previous 
considerations, where $P$ and $P^\vee$ were used
depending on the setting, 
{\em The indices of $X_b$ will be always from the 
lattice $P$ for both,
the twisted and untwisted settings}.

  
Instead of giving in this section the definition of $E$-dag 
polynomials and then restricting ourselves to their 
extremal parts, we will introduce the latter directly using 
mainly Proposition 2.5,$(i)$ from \cite{CO1}. Section 
\ref{SEC:GENPER} contains a systematic approach via
general $E$-polynomials.  

\subsection{\bf The \texorpdfstring{$\mathfrak{T}$}{T}-operator}
We set 
\begin{align}\label{whatxb}
X_{[b,\ze]}\!=\!q^\ze X_b,\  
wa(X_b)\!=\!X_{wa(b)}\!=\!q^{-(b,a)}X_{w(b)},
\ w\in W,\, a\in P,\!P^\vee,
\end{align}
where the latter formula generally may require fractional
powers of $q$. Given $b\in P_-$, we set
$M_b\equal\sum_{c\in W(b)} X_c$. 

For $0\le i\le n$ and any $b\in P$, let
\begin{align}
\T^\natural_i(X_b)\equal
\begin{cases}
s_i(X_b)+X_b, & \text{ if }(b,\al_i)<0,\\
X_b, & \text{ if }(b,\al_i)=0,\\
0, & \text{ if }(b,\al_i)>0.
\label{intmodsi}
\end{cases}
\end{align}
To establish a connection with \cite{CO1}, 
$$
\T^\natural_i=(\overline{T}_i')^\ast 
\hbox{ for (the extremal parts of) }
\overline{T}_i'
$$ given by (2.37) there and 
naturally extended to $i=0$. 
Note that the inequalities in (\ref{intmodsi}) are opposite to 
those in (2.37) due to applying ${}^\ast$, which sends
$X_b\mapsto X_b^{-1},\, q\mapsto q^{-1}$. 

Obviously, $\T^\natural_i(\T^\natural_i-1)=0$ for any $i$ because
$(s_i(b),\al_i)=-(b,\al_i)$. Recall that 
$([z,\ze],[z',\ze'])=(z,z')$, which  
is needed for $i=0$:
\begin{align}\label{intmodso}
\T^\natural_0(X_b)\equal
\begin{cases}
q^{(b,(\th^{\bep})^\vee)}\,s_{\th^{\bep}}(X_b)+
X_b, & \text{ if } (b,\th^{\bep})>0,\\
X_b, & \text{ if }(b,\th^{\bep})=0,\\
0, & \text{ if } (b,\th^{\bep})<0,
\end{cases}
\end{align}
where $\th^{\bep}=\th,\vth$ in the untwisted and twisted settings
correspondingly; recall that $\,\al_0=\al_0^\nu=[-\vth,1]$
in the twisted case and $\,\al_0=[-\th,1]$ in the untwisted case.

For any reduced
decomposition $\hat{w}=\pi_r s_{i_l}\cdots s_{i_1}$, where
$r\in O$ and $l=l(\hat{w})$, the product
$\T^\natural_{\hw}\equal \pi_r \T^\natural_{i_l}\cdots 
\T^\natural_{i_1}$
depends only on $\hw\in\, \hat{W},\hat{W}^\nu$ and not
on the particular choice of this decomposition.
This readily follows from \cite{CO1};
see Proposition 2.4 there and its proof.
We arrive at the following definition-theorem.

\begin{theorem}\label{EXTDAG}
Let $\check{\rho}^{\bep}=\pi_{\check{\rho}^{\bep}}
u_{\check{\rho}^{\bep}}$
be the decomposition from (\ref{xwo}) for
$\check{\rho}^{\bep}\equal\rho$ in the twisted ($\bep\,=
\nu$) and 
$\check{\rho}^{\bep}\equal\check{\rho}$ in the untwisted cases
($\bep\,=\varnothing$) correspondingly. 
Then $u_{\check{\rho}^{\bep}}=w_0$, and we set 

\begin{align*}
\pi_\rho^{\bep}\equal \pi_{\check{\rho}^{\bep}}\,=\,
w_0(\check{\rho}^{\bep})^{-1}\,=\,\check{\rho}^{\bep}\,w_0,\ \ \,
\mathfrak{T}^{\bep}\,\equal\, 
\T^{\natural}_{\pi_{\check{\rho}^{\bep}}},
\end{align*}
where $l(\pi_\rho^{\bep})=2(\rho,\check{\rho})-\#\{R_+\}$.
Then for $b\in P_-$,
\begin{align}\label{extrdag}
\E^\dag_b\equal
q^{(\check{\rho}^{\bep},b)}\mathfrak{T}^{\bep}(M_b)\ =\ 
\sum_{c\in W(b)}q^{-e_c}\,X_c \hbox{\, for proper \,}
e_c\in \Z_+,
\end{align}
where $e_c=0$ if and only if $c=b$.
We will call $\E^\dag_b$ extremal dag-polynomials and
$e(-b,w)\equal e_{w(b)}$ extremal dag-degrees.
\end{theorem}
{\it Proof.}
This theorem is actually a combination of Proposition 2.5$(i)$
and Corollary 2.6$(i)$ from \cite{CO1} (in the twisted setting).
The proof of the latter was not given there; let us provide it here.

The operators $\T^\natural_i$ do not change $M_b$ for $b\in P_+$
when $i>0$, as well as $\T^\natural_0$ and $\pi_r$ when $q=1$. 
Therefore $\mathfrak{T}^{\bep}_{q\mapsto 1}(M_b)=M_b$ and 
$\E^\dag_b(q\mapsto 1)=M_b$. 
for $b\in P_-$.

Applying $\mathfrak{T}^{\bep}$ to any monomial $X_a$ 
($a\in P$) will produce exactly the same monomials as for
$\mathfrak{T}^{\bep}_{q\mapsto 1}(X_a)$ but with some powers
$q^m$ for $m\le 0$ as their coefficients. This results directly
from formulas (\ref{intmodso}). Indeed, certain monomials can
be terminated in process of applying $\T^\natural_i$ due to 
the last line in (\ref{intmodsi}), but no other cancelations
can occur. It is the same process for $q\neq 1$ and for $q=1$;
only the resulting coefficients can be different.  
This argument can be equally used when calculating 
$\mathfrak{T}^{\bep}(M_b)$, and this justifies the decomposition 
in (\ref{extrdag}). The non-positivity of the $q$-degrees follows
automatically from the inequality in the first line
of (\ref{intmodso}).

Let us describe monic (with coefficient $1$)
monomials in $\E^\dag_b$ for $b\in P_-$.
Using the definition of $\T^\natural_i$,
one can always pick the term with $s_i$ in each
(\ref{intmodsi}) when applying $\mathfrak{T}^{\bep}$ to $X_{w_0(b)}$.
This can be seen directly, but it is more convenient to
switch there from $\T^\natural_i$ to the so-called $G$-operators
and use (2.6) from \cite{CO1}; see formula (\ref{rhogb}) below. 

We obtain that the following
monomial will be present in  $\mathfrak{T}^{\bep}(X_{w_0(b)})$:
\begin{align*}
\mathfrak{T}^{\bep}(X_{w_0(b)})\,\ni\, \pi_\rho^{\bep}(X_{w_0(b)})=
\check{\rho}^\bep(X_b)=q^{-(\check{\rho}^\bep,b)}X_b. 
\end{align*}
Due to the multiplier $q^{(\check{\rho}^\bep,b)}$
in (\ref{extrdag}), $X_b$ will have the coefficient $1$
in $\E^\dag_b$. And this is the only way to obtain monic
monomials in $\E^\dag_b$, which can be observed using the
same argument.

The fact that $X_b$ is a unique monomial in $\E^\dag_b$
with the coefficient that is not $q^m$ for $m<0$ can
be seen using the limit of the extremal part of
the nonsymmetric Macdonald polynomial $E_b(X;q,t)$ 
for $\,b\in P_-$ upon $\,t\to \infty$
and $\,q\to \infty$. Here one can
involve, for instance, the general theory of Matsumoto spherical 
functions, corresponding to the limit $\,q\mapsto \infty$ of
the $E$-polynomials (see \cite{Ch3}).
\sq
\medskip
 
\subsection{\bf Using \texorpdfstring{$G$}{G}-operators}
Let us apply the technique of $G$-operators, which
is standard in the theory of DAHA;
cf. (2.6) from \cite{CO1}. It helps to analyze $\E^\dag_b$ 
theoretically and is the best for practical calculations. For 
$\tal=[\al,j]$ and $b\in P$, let
\begin{align}
\G'_{\tal}(X_b)\equal
\begin{cases}
q^{j(b,\al^\vee)}\,s_\al(X_b)+X_b, & \text{ if }(b,\al)>0,\\
X_b, & \text{ if }(b,\al)=0,\\
0, & \text{ if }(b,\al)<0.
\label{intmodg}
\end{cases}
\end{align}

For a reduced decomposition 
$\pi^{\bep}_{\rho}=\pi_r s_{i_l}\cdots s_{i_1}$,
we set 
\begin{align}\label{alsuperseq}
\al^1=\al_{i_1},\, \al^2=s_{i_1}(\al_2),\, 
\al^3=s_{i_1}s_{i_2}(\al_3),\ldots, 
\al^l=-\pi^{\bep}_{\rho}\pi_r(\al_{i_l}).
\end{align}
Then the definition of $\E^\dag_b$
for $b\in P_-$ can be rewritten as follows:
\begin{align}\label{rhogb}
\E^\dag_b\ =\ q^{(\check{\rho}^\bep,b)}\,
\check{\rho}^\bep\, \G'_{\al^l}\,\G'_{\al^{l-1}}\cdots 
\G'_{\al^1}\,(M_b).
\end{align}
This formula readily gives that $\E^\dag_b$ contains the 
monomial $X_b$
(with the coefficient $1$) and it occurs only for $X_b$,
as it is stated in the theorem.
Indeed, we can always pick $X_b$
when applying (\ref{intmodg}) due to the fact that the nonaffine
parts of all $\al^j$ are negative (so their affine extensions 
are strictly positive integers). Then 
$q^{(\check{\rho}^\bep,b)}\check{\rho}^\bep(X_b)=X_b$ and
there will be always nontrivial powers of $q$ for any
other choices.

Let us check (\ref{rhogb}). Generally, the $\,G$-operators\,
extend $\{G_{\al_i}=s_iT_i, \, 0\le i\le n\}$ from $\{\al_i\}$
to arbitrary $\tal\in \tR^\bep$. Here 
the obvious relations $s_i \T^\natural_i=\T^\natural_i\,$
hold, so we can we extend $\T^\natural_i$ themselves:
\begin{align}
\T_{\tal}^\#(X_b)\equal
\begin{cases}
s_{\tal}(X_b)+X_b, & \text{ if }(b,\al)<0,\\
X_b, & \text{ if }(b,\al)=0,\\
0, & \text{ if }(b,\al)>0.
\label{intmodgti}
\end{cases}
\end{align}
Thus $\T^\natural_{\tal}$ are
given by formulas (\ref{intmodg})
upon $q\mapsto q^{-1}$ and for inequalities from (\ref{intmodsi}),
i.e. opposite to those in (\ref{intmodg}). The following is 
straightforward:
\begin{align}\label{edagbt}
\E^\dag_b\ =\ q^{(\check{\rho}^{\bep},b)}
\check{\rho}^\bep w_0 \,\T^{\natural}_{\al^l}\,
\T^\natural_{\al^{l-1}}\cdots 
\T^\natural_{\al^1}\,(X_{w_0(b)}).
\end{align}
Then we move $w_0$ to the right using that
$w_0(\check{\rho}^\bep)=-\check{\rho}^\bep$,
which implies that the automorphism $-w_0$ transforms  
$\pi^{\bep}_{\rho}=\pi_r s_{i_l}\cdots s_{i_1}$ to another reduced
decomposition of $\pi^{\bep}_{\rho}.$ The product of $\T^\natural$
in (\ref{edagbt}) does not depend on the choice of the reduced 
decomposition; thus (\ref{rhogb}) is checked.
\smallskip

\subsection{\bf Total additivity}
The following is the key in establishing the connection to
the PBW-filtration.

\begin{theorem}\label{EXTADDIT}
For arbitrary $b,c\in P_+$ (we say, totally), 
the additivity of the extremal dag-degrees holds\,:
\begin{align*}
e^{\bep}(b + c,w)=e^{\bep}(b,w)+e^{\bep}(c,w),
\hbox{\,\ where\, } w\in W
\end{align*}
and (as above) backepsilon $\bep\,$
means $\,\nu\,$ in the twisted case and
$\varnothing$ in the untwisted case.
\end{theorem}
{\it Proof.} 
We will expand the argument that provided the
pure $q$-powers in the expansion from
(\ref{extrdag}). The process of applying consecutive
$\T^\natural_i$ when calculating $\mathfrak{T}^{\bep}(M_b)$ 
for $b\in P_-$
is either by adding $s_{\al_i}(q^m X_d)$ for any existing 
monomial $q^m X_d$, or annihilating it or leaving it
unchanged. 
Using  the passage to $q=1$, as when establishing (\ref{extrdag}), 
we obtain that for each $a \in W(b)$,
there exists a {\em unique $a_\circ\in W(b)$} and 
a {\em unique
sequence of selections of either 
$s_{i_p}(q^m X_d)$ or $q^m X_d$\,} from the first line
of (\ref{intmodsi}) at each  $\,s_{i_p}\,$ 
in $\,\pi^{\bep}_\rho=\pi_rs_{i_l}\cdots s_{i_1}\,$ satisfying 
$(\al_{i_p},d)<0$
such that the resulting monomial from 
$\mathfrak{T}^{\bep}(X_{a_\circ})$ is nonzero and proportional 
to $X_a$.
Recall that the places where $(\al_{i_p},d)\ge 0$ do not change 
$q^m X_d$ or annihilate it. 

The inequalities $(d,\al_{i_p})< 0$ from
(\ref{intmodsi}) and $i=i_p\,$
can be recalculated
to the form $(b,\be)<0$ for proper $\be\in R_+$ (depending on 
the particular step and the previous selections). This root
$\be$ must be positive. Indeed, if it is not, then $(b,\be)=0$
since $b\in P_-$. Note that if $(b,\be)=0$, then
we can pick an arbitrary term (from two) in the first line 
of (\ref{intmodsi}) in this case (they coincide), so the 
procedure is uniform for  $(b,\be)\le 0$. 

Assuming now that $(b,\al)<0$ for all $\al\in R_+$,
we see that this sequence of selections
is uniquely determined by $u,w\in W$ such that 
$u(b)=a_\circ, w(b)=a$. We conclude that if the resulting 
monomial is nonzero for one $b$, then it is nonzero for all 
$b\in P_-$ including those from the boundary of the negative 
Weyl chamber.  

Finally, the resulting coefficient of 
$X_a$ is the product of the (negative) powers of $q$ calculated
at all elements $s_0$ in the reduced decomposition of 
$\pi^{\bep}_\rho$ where the term $s_0(q^m X_d)$ was selected. 
Obviously, its $q$-degree is a linear function of $b\in P_-$, 
which proves the required additivity. This leads to 
an algorithm of finding $q$-degrees, which we hope to discuss
in further works.
\sq
\smallskip

Proposition 2.5 from \cite{CO1} can be extended to prove
the following stronger version of Theorem \ref{EXTADDIT}
based on the exact analysis of the elements $a_\circ=u(b)$
that appeared in its proof. 

\begin{theorem}\label{ADDWITHIN}
For $b\in P_+$,
the polynomial $\,\mathfrak{T}^{\bep}(X_{u(b)})\,$ is 
nonzero if and only if $\,u=$id\, or the product
of pairwise commutative simple (nonaffine) reflections.
Let $\,e^{\bep}(b,w;u)\,$ be $\,e^{\bep}(b,w)\,$ 
if the monomial $\,q^{m}X_{w(b)}\,$ for $m=e^{\bep}(b,w)$
has a nonzero coefficient in
$\,\mathfrak{T}^{\bep}(X_{u(b)})\,$  and
zero otherwise. Then $e^{\bep}(b,w;u)\neq 0\,$ occurs exactly
for one such $u$ modulo $\,W^b\,$ and one has: 
\begin{align*}
e^{\bep}(b + c,w;u)=e^{\bep}(b,w;u)+e^{\bep}(c,w;u),
\hbox{\,\ where\, } w\in W,\ b,c\in P_+. 
\end{align*}
\sq
\end{theorem}

It is not too difficult to obtain the formulas for $e(b,w)$
for anti-fundamental weights $b=-\om_i$ for classical root systems; 
their calculation is based on a relatively straightforward induction 
with the respect to $l(w)$ (to be continued in further works).
The formulas for the exceptional root systems were calculated 
mainly using computers; they are long for $E_{7,8}$. We provide
them for $G_2$, $F_4$ and (partially) for $E_6$. 
\medskip

\subsection{\bf The case of \texorpdfstring{$A_n$}{An}}
All fundamental weights are minuscule for
$A_n$ ($\mathfrak{g}=\msl_{n+1}$) and we do not actually
need the $\mathfrak{T}$-operator  to calculate $\E^\dag_{-\om_i}$.
One can directly send $t\to\infty$ in following special case of the 
Haiman-Haglund-Loehr formula
\cite{HHL} for $E$-polynomials.
The variables $x_i\, (1\le i\le n+1)$ correspond 
to $\vep_i$ from the $A_n$-table from \cite{Bo}. 
\begin{proposition} 
\begin{align}\label{glnepol}
&E_{-\om_i}\ =\sum_J\,
 x_{\!J}\,\prod_{k=1}^{n_i(J)}\,\frac{1-t^k}{1-qt^k} 
\hbox{\ \, for \ \,} \#\{J\}=i,\\
J=&\{1\!\le\! j_1\!<,\ldots,<\!j_i\!\le\! n+1\},\ 
n_i(J)\!=\!\#\bigl\{J\cap [1,n+1-i]\bigr\},\notag
\end{align}
where $\#$ is the cardinality of a set and 
$x_{\!J}=\prod_{k=1}^i x_{j_k}$.
\end{proposition}
{\it Proof.}
To prove the formula \eqref{glnepol} we fix a number 
$k\le \min(i,n+1-i)$ and look at the coefficient 
in $E_{-\om_i}$ of the monomial
\[
x_{\!J}=x_{n+1}\dots x_{n-k+2}\,x_1\dots x_{i-k}  
\]
(we note that $n_i(J)=k$, since $n-k+2\ge i+1$).
We use Theorem 3.5.1 (formula (26)) of \cite{HHL} to compute
this coefficient; the notation is from this paper. The 
{\em composition} $\mu$ in our case is simply 
$(0,\dots,0,1,\dots,1)$ with $i$ units. Let us denote the
entries of a {\em filling\,} $\sigma$ of $\mu$ by 
$\sigma_i,\dots,\sigma_1$ from left to right. We want to
find all {\em non-attacking\,} $\sigma$ such that $x^\sigma=x_{\!J}$.
They are given by
\[
\sigma_1=n+1,\dots,\sigma_k=n-k+2 \text{\, and \,} 
\{\sigma_{k+1},\dots,\sigma_i\}=\{1,\dots,i-k\}.
\]
Therefore the summation in the (special case of the) HHL formula 
runs  over the permutation group $\S_{i-k}$; for a permutation $g$,
we set $\sigma_{k+i}=g(i)$. For any $g\in \S_{i-k}$ and the 
corresponding $\sigma$, one has $\rm{maj}(\hat\sigma)=0$ and 
$\,\rm{coinv}(\hat\sigma)\,$ is the number of inversions in $g$
(here $\hat\sigma$ is the {\em augmented filling}).
The factor 
\[
\prod\frac{1-t}{1-q^{l(u)}t^{a(u)}}
\]
in their paper equals  
$(1-t)^{i-k}/\bigl((1-qt)(1-q^2t)\dots(1-q^{i-k}t)\bigr)$
for any such $g$. 
We thus obtain that the coefficient of $x_{\!J}$ in  
$E_{-\om_i}$ equals 
\[
\prod_{a=1}^{i-k-1} (1+t+\dots+t^a)\prod_{a=1}^{i-k}
\frac{1-t}{1-q^at}=
\prod_{a=1}^{i-k}\frac{(1-t^a)}{(1-q^at)}.
\]

To complete the proof of the proposition, we recall the invariance of 
$E_{-\om_i}$ with respect to the action of the product of symmetric
groups $\S_i\times \S_{n-i+1}$. 
\sq
\medskip

As an immediate application,
\begin{align}\label{glnepolinf}
\E^\dag_{b}=E^\dag_{b}=E_{b}\mid_{t\to\infty}\, = \, 
\sum_J\, q^{n_i(J)}x_{\!J}\, =\, \sum_{c\in W(b)} q^{n(-b,w)}X_c,
\end{align}
where $b=-\om_i,\, c=w(b\,)$ and $n(-b,w)$ depends only on $b-c$
(see the next subsection). 
The calculation of $\E^\dag_{-\om_i}$ for $A_n$ 
is simple to perform using directly (\ref{extrdag}) and 
this approach can be extended to any classical root systems.
We will discuss a systematic combinatorial theory of
the operator $\mathfrak{T}^{\bep}$ and the calculations
for the fundamental weights elsewhere.

\comment{
For $A_n$, one can send $t\to\infty$ in
following special case of the Haiman-Haglund-Loehr formula
\cite{HHL} for $E$-polynomials.
Let us use the variables $Z_i\, (1\le i\le n+1)$ corresponding 
to $\vep_i$ from the $A_n$-table from \cite{Bo}. 
Setting $Z_1\cdots Z_{n+1}=Z_\circ$ and
$[i+1,n+1]=\{i+1,\ldots,n+1\}$,
\begin{align}\label{glnepol}
&E_{-\om_i}\ =\ Z_\circ^{i/{n+1}}\, \sum_J\,
\prod_{k=1}^{n_i(J)}\,\frac{1-t^k}{1-qt^k}\ Z_J/Z_\circ 
\hbox{\ \, for } \\
J=&\{1\!\le\! j_1\!<,\ldots,<\!j_{n-i+1}\!\le\! n+1\},\ 
n_i(J)\!=\!\#\bigl\{J\ominus [i+1,n+1]\bigr\},\notag
\end{align}
where $\#$ is the cardinality of a set, $\ominus$ the symmetric
difference of two sets and $Z_J=\prod_{k=1}^{n-1+1}Z_{j_k}$.
For $J=[i+1,n+1]$, the coefficient of 
$Z_\circ^{i/{n+1}}Z_{J}/Z_\circ=X_{-\om_i}$ in $E_{-\om_i}$
is $1$, and it is the only monic monomial in this sum.
Generally, $n_i(J)$ equals the minimal length of 
the decomposition of  $\om_i-w(\om_i)$ in terms of positive 
roots for $J=w([i+1,n+1])$ and any $w\in \S_{n+1}$. It will
be denoted by $n(\om_i,w)$ in the next section. Indeed,
$Z_J/Z_{[i+1,n+1]}$ is the product of $n_i(J)$ monomials 
$X_\al=Z_l/Z_m$ for $\al=\vep_l-\vep_m$, where $l<m$
and all numbers $\{l,m\}$ are pairwise different.
Finally,  for $b=-\om_i$, where $1\le i\le n$, one has
\begin{align}\label{glnepolinf}
&\E^\dag_{b}=E_{b}\mid_{t\to\infty}\,=\, 
Z_\circ^{i/{n+1}}\sum_J\,
q^{n_i(J)}Z_J/Z_\circ\, =\, \sum_{c\in W(b)} q^{n(-b,w)}X_c,
\end{align}
where $c=w(b)$ and $n(-b,w)$ depends only on $b-c$. The
calculation of $E_{-\om_i}$ can be performed using 
directly (\ref{extrdag}) in this and any classical cases.
We will discuss the details elsewhere.
}
\medskip

\setcounter{equation}{0}
\section{\sc Kostant \texorpdfstring{$q$}{q}-partition function}
In this and the next section we switch from using $\,b,c\in P$
to the standard $\la,\mu\in P$ in the Lie theory. The main 
reason of this split of notation is that $b$ is mainly antidominant
in the theory of Macdonald polynomials, which corresponds dominant
$\la$ in what will follow. Also, we used $\,b$ to ensure the
maximal compatibility with \cite{CO1} (and quite a few other papers 
on the Macdonald polynomials).

The definition of the {\em extremal $q$-degrees}
of Lusztig's $q$-analogous of Kostant partition
function is as follows.
Let $n(\la,w)$ for $\la\in P_+$ and $w\in W$
be the minimal number of terms
in the decomposition of $\la-w(\la)$ in terms
of positive roots. In the twisted setting,
we count long roots with multiplicity $\nu_{\lng}$.
Recall that $\nu_\al=(\al,\al)/2$ and $\nu_{\sht}=1$.
To avoid possible confusions we will frequently (but 
not always) use the notation $n^\nu(\la, w)$ in the twisted case. 
Let us begin with considering some simple examples.

\subsection{\bf The case of reflections}
In the untwisted setting, one has 
\begin{align}\label{nbsunt}
n(\la, s_\al)=(\la,\al^\vee) \hbox{ \, for any \, } 
\al\in R_+,\, \la\in P_+,
\end{align}
except for $G_2$ and
$\al=\al_1+\al_2$ or $\al=2\al_1+\al_2$
(in the notation from \cite{Bo}).
Let us check this. Since 
$\la-s_\al(\la)=m\al\,$ for $\,m=(\la,\al^\vee)$, 
$\,n(\la,s_\al)\neq (\la,\al^\vee)\,$ if and only if 
$m \al=\sum_{j=1}^{M}\be^j$ for $M\le m-1$ and
certain positive roots $\be^j$. Reducing $\,m\,$ if
necessary in this relation, we can assume 
that $\be^j\neq \al$. Then 
\begin{align}\label{malbe}
(m\al,\al^\vee)=2m=\sum_{j=1}^{M}
(\be^j,\al^\vee)\,\le\, \sum_{j=1}^{M}\nu_{\be_i}/\nu_\al.
\end{align}
Due to the last inequality,
(\ref{malbe}) can be valid only for short $\al$; moreover,
all $\{\be^j\}$ must be long and also the root system
must be $G_2$. 
For $G_2$, $\,n(\la,s_\al)\neq (\la,\al^\vee)\,$ occurs 
due to the following relations: 
$$
3(\al_1+\al_2)\!=\!(3\al_1+2\al_2)+\al_2 \hbox{\, or\, } 
3(2\al_1+\al_2)\!=\!(3\al_1+2\al_2)+(3\al_1+\al_2).
$$
Hence the exceptional cases are $n(\la,s_\al)=2k+r$ when 
$\al=\al_1+\al_2$ or $\al=2\al_1+\al_2$ and 
$(\la,\al^\vee)=3k+r,\, 0\le r\le 2$.
Otherwise (\ref{nbsunt}) holds.
\smallskip

{\em The twisted setting}. Now long roots are counted with
multiplicity $\nu_{\lng}$ and (\ref{malbe})
readily gives that $\,n^\nu(\la,s_\al)=(\la,\al^\vee)\,$ 
for any short $\al>0$.

For long $\al>0$, let $k$ be the number of long roots $\be^j>0$
in the minimal possible presentation $m \al=\sum_{j=1}^{M}\be^j$.
If  $\,n^\nu(\la,s_\al)\neq \nu_{\lng}(\la,\al^\vee)\,$, then  
$M$ must be smaller than $(m-k)\nu_{\lng}$. As above, $m\,$ is
$\,(\la,\al^\vee)\,$ reduced by the number of $\be^j=\al$; so
we assume that $\be^j\neq \al$.
We arrive at the inequality 
\begin{align}\label{malbet}
(m\al,\al^\vee)=2m=\sum_{j=1}^{M}
(\be^j,\al^\vee)\,\le\, M,
\end{align}
which can hold (again) only in the case
of $G_2$. For $G_2$, the exceptional cases are 
$\,n^\nu(\la,s_\al)=2(\la,\al^\vee)<\nu_{\lng}(\la,\al^\vee)\,$ 
for long $\al\neq \al_2$. 

Finally, except for $G_2$, the formula for $\,n^\nu\,$ 
reads as 
\begin{align}\label{nnubst} 
n^\nu(\la,s_\al)=\nu_\al(\la,\al^\vee),\ \la\in P_+, \al\in R_+.
\end{align}
\smallskip

\subsection{\bf Maximal roots}\label{maxroots}
Let us consider now $\la=\th,\vth$ (they are dominant).
Then $n(\th,w)=1$ and $n^\nu(\th,w)=\nu_{\lng}$ 
for the maximal long root $\th$, correspondingly,
in the untwisted and twisted setting 
provided $(w(\th),\th)>0$ and $w(\th)\neq\th$
and excluding the case $\th-w(\th)=2\al_1+\al_2$ for 
twisted $G_2$. Indeed, then $\th-w(\th)$ is a long (positive) 
root due to $(w(\th),\th)>0$ unless $w(\th)=\th$.
When $\th-w(\th)=2\al_1+\al_2$ for twisted $G_2$, one has 
$2\al_1+\al_2=\al_1+(\al_1+\al_2)$ and $n^\nu(\th,w)=2$ in this case.

Switching here to the short maximal root $\vth$
and imposing the conditions $(w(\vth),\vth)>0$ and   
$w(\vth)\neq\vth$ (otherwise, $n(\vth,w)=0$), one obtains  
that $n(\vth,w)=1=n^\nu(\vth,w)$ for any, twisted or
untwisted, setting (including $G_2$). Indeed,
$\vth-w(\vth)$ is a short root in this case.
\smallskip

Let us impose now the opposite inequality 
$(w(\th),\th)<0$ and check that
$n(\th,w)=2$. In this case,
$\al=-w(\th)\in R_+$ and $\th-w(\th)=\th+\al$ is a 
sum of two positive (long) roots; the latter
can not be a (single)  root due to the maximality of $\th$. 
If $\nu_{\lng}\neq 1$ (i.e. $R\,$ is not simply-laced),
then $\th+\al$ cannot be a sum $\be+\ga$ of two positive roots 
where $\be$ (one of them) is short, since otherwise  
$
(\th+\al,\th+\al)=8\nu_{\lng}=(\be+\ga)^2\le 4+6\nu_{\lng}.
$
This gives that  $n(\th,w)=2$ in the untwisted case.

Continuing this argument, $\th+\al$ 
cannot be a sum of $3$ positive short roots, 
since otherwise $|\th+ \al|^2\le 12$. This gives
that  $(w(\th),\th)<0$ results in
$n^\nu(\th,w)=2\nu_{\lng}$ in the twisted case 
except for $G_2$.

Similarly,
the condition $(w(\vth),\vth)<0$ 
results in $n(\vth,w)=2=n^\nu(\vth,w)$ for both 
settings and including the root system $G_2$.
Indeed, if the difference $\vth-w(\vth)$ is a single root
then it must be short, which
contradicts the maximality of $\vth$ among short roots. 
\smallskip

The remaining cases are when $(w(\th),\th)=0$ or
$(w(\vth),\vth)=0$. Let us check that
$n(\th,w)=2$ under the first equality
in the untwisted setting. Indeed, 
$\th-w(\th)$ is not a root since any
sum/difference of two pairwise orthogonal long
roots can not be a root; so  $n(\th,w)\ge 2$. 
It is obviously exactly $2$ if $w(\th)<0$. 

If $w(\th)>0$, then one can find
$\be\in R_+$ such that $w(\th)+\be\in R_+$ and 
$(w(\th)+\be,\th)>0$. 
This gives that
$\be'=w(\th)+\be-\th$ is a root from $R_-$,  
$\th-w(\th)=\be-\be'$ and $n(\th,w)=2$ in this case.

We use here that the condition
$(w(\th),\th)=0$ simply means that supp$(w(\th))$, 
a connected set in $\Ga$ formed by the simple roots that
occur in the expansion of $w(\th)$, 
does not contain the simple roots (one or two)
adjacent to $\al_0\in\tilde{\Ga}$. 
Then we connect supp$(w(\th))$
with $\al_0$ by a segment; the sum of the
simple roots in this segment (excluding $\al_0$ 
and supp$(w(\th))$) gives $\be$.
\smallskip

Similarly, $(w(\vth),\vth)=0$ implies $n^\nu(\vth,w)=2$
in the twisted case. First of all, this condition is 
empty for $G_2$. 
Then $\vth=\vep_1$ for $B_n(n\ge 2)$ and $F_4$
in the notation from \cite{Bo}; therefore $\vth-w(\vth)$ is
always a (single) long root or $2\vep_1$ in these cases 
(due to $(w(\vth),\vth)=0$). Finally, $\vth=\vep_1+\vep_2$ for
$C_n(n\ge 2)$ and $(w(\vth),\vth)=0$ if and only if
$w(\vth)=\pm(\vep_1-\vep_2)$ or $w(\vth)=\pm\vep_i\pm\vep_j$
for $i,j>2, i\neq j$. Thus $\vth-w(\vth)$ is a 
single long root or a sum of two short positive roots
for $C_n$ as well as for $B,F$.

We leave the consideration of $\la=\th$ in the twisted setting
and $\la=\vth$ in the untwisted setting to the reader
(correspondingly 
under $(w(\th),\th)=0$ and $(w(\vth),\vth)=0$).
We arrive at the following lemma.

\begin{lemma}\label{THETAW}
Let $\th'$ be $\th$ or $\vth$; we will exclude the
case of $G_2$ if $\th'=\th$ in the twisted setting.
Then
\begin{align}
&n(\th',w)=1 \hbox{\, and \,} n^\nu(\th',w)=\nu_{\th'}
\hbox{\, if\, } (w(\th'),\th')>0 \hbox{ and \,} 
w(\th')\neq \th',\notag\\
&n(\th',w)=2 \hbox{\, and \,} n^\nu(\th',w)=2\nu_{\th'}
\hbox{\ \ assuming that\ \ }\, (w(\th'),\th')\,<\,0,\notag \\
&n(\th,w)=\,2\,= n^\nu(\vth,w) \hbox{\,\, if\, } 
(w(\th),\th)= 0=(w(\vth),\vth)
\hbox{\, correspondingly.}\notag
\end{align}
\sq
\end{lemma}

Similar direct analysis can be used for minuscule weights,
which we will omit in this paper. This is directly
related to the fact that PBW- and dag-degrees coincide with
$n(\om,w)$ for minuscule $\om$ and when $\om=\th'$ 
(sometimes even for both, $\th$ and $\vth$), which follows
from  the formulas we provide below.
\smallskip

\comment{
Let $w'=s_iw$ provided $l(w')=l(w)+1$. Then
\begin{align*}
\la-w'(\la)=\la-w(\la)+(w(\la),\al_i^\vee)\al_i.
\end{align*}
Here $w^{-1}(\al_i)\in \la(s_i w)$;
for instance $w^{-1}(\al_i)$ is positive
and $(w(\la),\al_i^\vee)=(\la,w^{-1}(\al_i^\vee))
\ge 0$. Thus we obtain that
$n(\la,w')-n(\la,w)\le (\la,w^{-1}(\al_i^\vee))$.
The inequality is strict if $(\la-w(\la),\al_i)<0$.
Indeed,  $(\be_j,\al_i)<0$ in this case at least 
for one $\be_j$ in any {\em minimal} possible decomposition
$\la-w(\la)=\sum_j \be_j$  with positive $\be_j$
and $1\le j\le n(\la,w)$. Thus $\be_j+\al_i\in R_+$
for such $j$. Moreover, it is strict if  $(\la-w(\la),\al_i)=0$
but $(\al_i,\be_j)\neq 0$ for at least one $\be_j$
from any minimal decomposition of $\la-w(\la)$.
This approach can be used, for instance, to determine 
$n(\la,w)$ for minuscule weights $b$.
}

\subsection{\bf Extremal additivity}
\begin{theorem}\label{ABDaddit}
The additivity $n(\la,w)+n(\mu,w)=n(\la+\mu,w)$ holds for
arbitrary  $\la,\mu\in P_+$ and any $w\in W$ for the
following root systems:
$$
A_n (n\ge 1), \hbox{\, untwisted } C_n (n\ge 2),
\hbox{\, twisted } B_n (n\ge 2).
$$
Moreover, the total additivity (any $w\in W, \la\in P_+$) 
holds only for these root systems and twisted $G_2$.
\end{theorem}
{\it Proof.} 
The counterexamples for the total additivity (the second part 
of the theorem) will be given below. Let us prove  
the first part.
\smallskip

We start with the case of $A_n$.
Then the positive roots are
\[
\al_{i,j}=\al_i+\al_{i+1}+ \dots +\al_j, \ 1\le i\le j\le n,
\]
where $\al_{i,i}=\al_{i}$ are the simple roots. The Weyl group 
is equal to the symmetric group $\S_{n+1}$. The fundamental weights
are denoted by $\omega_i$, $i=1,\dots,n$.
 
Let $\la=\sum_{i=1}^n m_i\omega_i$. It is convenient to pass from 
the $\msl_n$-weights to the $\mgl_n$-weights. To this end, we define 
$\la_i=m_1+\dots +m_{n+1-i}$, $1\le i\le n$.
Then we have $\la_1\ge \la_2\ge \dots\ge\la_n$ and this is exactly 
the $\mgl_n$-weight of the highest weight vector. We write 
\begin{equation}\label{dif}
\la-w(\la)=(\la_1-\la_{w(1)},\la_2-\la_{w(2)},\dots,\la_n-\la_{w(n)}).
\end{equation}
We note that $\al_{ij}=\varepsilon_i-\varepsilon_j$ and therefore our 
task is as follows. We must write the right hand side of \eqref{dif} 
as a sum of $n$-tuples, corresponding to $\al_{i,j}$, 
i.e. with $1$ at the $i$-th place, $-1$ at the $j$-th place and 
zeros elsewhere, minimizing the number 
$s$ of summands. Obviously, $s$ equals the sum 
of nonnegative terms on the right-hand side of 
\eqref{dif}, which can be readily calculated: 
\[
\sum_{i:\  \la_i\ge \la_{w(i)}} (\la_i-\la_{w(i)}) = 
\sum_{i:\  w(i)>i} (m_i +\dots + m_{w(i)-1}).
\]     
Clearly, this expression is linear in $m_i$'s, 
which gives the required.
\smallskip

{\em The $C_n$-case.} 
Let us fix pairwise orthogonal weights $\veps_1,\dots,\veps_n$. 
The positive roots of $\msp_{2n}$ are of
the form $\veps_i-\veps_j$, $i<j$ and $\veps_i+\veps_j$, $i\le j$. 
The fundamental weights are given by
\[
\omega_i=\veps_1+\dots+\veps_i,\ i=1,\dots,n.
\]
Hence any weight $\la=m_1\omega_1+\dots+m_n\omega_n$ can be presented 
as a Young diagram 
\[
\la=(\la_1\ge\dots\ge\la_n),\ \la_i=m_i+\dots+m_n.
\]
The Weyl group $W$ contains all the permutations from $\S_n$ as 
well as all sign changes $\veps_i\to -\veps_i$.
So we can represent each element of $w\in W$ as a map from the set 
$\{1,\dots,n\}$ to the set $\{1,\dots,n,-n,\dots,-1\}$.

We want to prove that 
\begin{multline*}
n(\la,w) = \sum_{w(i)>i} (m_i+\dots +m_{w(i)-1}) + 
\sum_{w(i)<0} (m_i+\dots + m_n)\\
=\sum_{w(i)>i} (\la_i -\la_{w(i)}) + \sum_{w(i)<0} \la_i.
\end{multline*}
This formula actually follows from that in the $A_{2n-1}$-case. 
Namely, to any dominant $\msp_{2n}$-weight
$\la=(\la_1,\dots,\la_n)$, we associate the $\msl_{2n}$-weight 
$\tilde\la$ defined by adding $n$ zeros after $\la_n$.
Also, given $w\in W$ (the Weyl group of type $C_n$), 
we associate with it $\tilde w$ from the Weyl group for 
$\msl_{2n}$ defined as follows.
If $w(k)>0$, then $\tilde w(k)=w(k)$. If $w(k)<0$, then 
$\tilde w(k)=2n+1-k$. All other values of $\tilde w$ are
insignificant.  Now assume that $\la-w(\la)$ is decomposed as 
a sum of positive roots of $\msp_{2n}$ and the number of summands 
is the minimal one. Obviously the roots 
\[
\veps_i-\veps_j \text{\,\, and\,\, } \veps_k+\veps_j
\] 
can not appear in this decomposition simultaneously; 
otherwise they can be summed up to a single root. Therefore each
$\veps_j$ enters such decomposition (of minimal possible
length) with the same
sign (in the corresponding positive roots) and no cancelations 
occur.

Now let us attach to the $\msp_{2n}$-roots in the form 
$\veps_i+\veps_j$ 
the $\msl_{2n}$-roots $\veps_i - \veps_{2n+1-j}$ and to the 
$\msp_{2n}$-roots in the form $\veps_i-\veps_j$ the $\msl_{2n}$-roots
$\veps_i - \veps_{j}$. Then the minimal length decomposition of 
$\la-w(\la)$ into a sum of $\msp_{2n}$ positive roots
induces the decomposition of $\tilde\la-\tilde{w}(\tilde \la)$ into a 
sum of positive $\msl_{2n}$-roots.

In the opposite direction, a decomposition of 
$\tilde\la-\tilde{w}(\tilde \la)$ induces a decomposition of
$\la-w(\la)$. Hence   
\[
d(\la,w) = d(\tilde \la,\tilde w)=
\sum_{w(i)>i} (\la_i -\la_{w(i)}) + \sum_{w(i)<0} \la_i. 
\] 
Obviously, this expression is linear in terms of $\la$. 
\smallskip

The twisted $B_n$-case is similar to the untwisted 
$C_n$-case. The twisted $G_2$ is actually similar
to $A_2$; the shorts roots mainly
appear in the minimal decompositions because long ones
are counted with the multiplicity $\nu_{\lng}=3$. We will publish 
the details elsewhere.   


\subsection{\bf Counterexamples to additivity}\label{Cta}
Addressing the second part of the theorem,
let us provide the examples when 
\begin{align}\label{bwball}
n(\la,w)+n(\mu,w)\neq n(\la+\mu,w)\hbox{\, for some\, } 
w\in W,\, \la,\mu\in P_+;
\end{align}
then this inequality
can be only in the following direction: 
$n(\la,w)+n(\mu,w)> n(\la+\mu,w)$.

{\em For untwisted $B_3$} in the notation from \cite{Bo},
one can take $\la=\om_1+\om_3$ and $w=w_0$. 
Abbreviating $A[a,b,c]=A[abc]\equal a\al_1+b\al_2+c\al_3$,
\begin{align}\label{bwbb3}
&\la-w(\la)=A[111]+A[112]+A[122],\ n(\la,w)=3,\\
&\om_1-w(\om_1)=A[100]+A[122],\, n(\om_1,w)=2,\notag\\ 
&\om_3-w(\om_3)=A[001]+A[122],\, n(\om_3,w)=2. \notag
\end{align}
Here all $A[\ldots]$ are positive roots and therefore
$n(\la,w)<n(\om_1,w)+n(\om_3,w)$. Using the standard
embeddings, this provides counterexamples for all 
untwisted $B_n(n>3)$ and untwisted $F_4$.
\smallskip

{\em For twisted $C_3$} (i.e. that with $\,n^\nu\,$),
let $\la=\om_1+\om_3$ and $w=w_0$. Then
\begin{align}\label{bwbc3}
&\la-w(\la)\,=\,2A[111]+2A[121],\ \ n^\nu(\la,w)=4,\\
&\om_1-w(\om_1)=A[100]+A[121],\, n^\nu(\om_1,w)=2,\notag\\ 
&\om_3-w(\om_3)=A[011]+A[111]+A[121],\, n^\nu(\om_3,w)=3,\notag 
\end{align}
where $A[\ldots]$ are all positive short roots and 
$n^\nu(\la,w)<n^\nu(\om_1,w)+n^\nu(\om_3,w)$.
This automatically provides examples of (\ref{bwball}) for all 
twisted $C_n(n\ge 3)$ and twisted $F_4$. Note that
$\om_1-w(\om_1)=A[221]$ is a (single) long root, so it is counted
as $1$ in the untwisted setting; so the equality 
$n(\la,w)=n(\om_1,w)+n(\om_3,w)$ holds in the untwisted case
for this $\la$.
\smallskip

{\em In the case of $D_4$},
let $\la=\om_3+\om_4,\, w=w_0.$ Then
\begin{align}\label{bwbd4}
&\la-w(\la)=A[0111]+A[1111]+A[1211],\ n(\la,w)=3,\\
&\om_3-w(\om_3)=A[0010]+A[1211],\, n(\om_1,w)=2,\notag\\ 
&\om_4-w(\om_4)=A[0001]+A[1211],\, n(\om_3,w)=2, \notag
\end{align}
where $A[\ldots.]$ are positive roots. This gives examples
of (\ref{bwball}) for any $D_n(n\ge 4)$ and
$E_{6,7,8}$.
\smallskip

{\em For untwisted $G_2$}, let $\la=c_1\om_1+c_2\om_2$.
The simplest weight when
\begin{align}\label{bwbg2}
n(\la,w)<c_1 n(\om_1,w)+c_2 n(\om_2,w),\ \,c_1,c_2\ge 0,
\end{align}
is $\la=2\om_1+\om_2$. The corresponding $w$ are
$w_0$ or $s_{2\al_1+\al_2}$. For $w=w_0$,
$\la-w(\la)=14\al_1+8\al_2=A[21]+A[31]+3A[3,2]$, which
gives $n(\la,w)=5$. However,
$\om_1-w(\om_1)=A[10]+A[32]$ and
$\om_2-w(\om_2)=2A[32]$, 
which makes the right-hand side of (\ref{bwbg2}) equal
to $6$. For $w=s_{2\al_1+\al_2}$, one has
\begin{align}\label{bwbg21}
&\la-w(\la)=14\al_1+7\al_2=A[21]+2A[31]+2A[3,2],\\
&\om_1-w(\om_1)=A[10]+A[32],\,
\om_2-w(\om_2)=A[31]+A[32]. \notag
\end{align}

\smallskip

\subsection{\bf Fundamental weights}
It is not difficult to calculate $\,n(\la,w)\,$ for
fundamental weights $\la=\om_i\, (1\le i\le n)$
and the corresponding minimal decompositions of $\,\om_i-w(\om_i)$ 
for the classical root systems (for any $w\in W$). 

We represent $w\in \S_{n+1}$ as
permutations $w=(w_j, 1\le j\le n+1)$ and
use the permutations with signs
$w=(w_1,\cdots, w_n)$ for $w\in W$ in types $BC\!D$. 
Recall that $(|w_i|:\, 1\le i\le n)$ is a permutation
of $\{1,2,\cdots,n\}$; the signs of $w_i$ 
can be arbitrary for $B_n,C_n$ and with
even number of minuses for $D_n$. 

We use $\#\{\,.\,\}$ for 
the number of elements of a given set and $\bigl[\,.\,\bigr]$ for 
the integer part. The formulas below for $a_i(w)$ and
(later) $\tilde{a}_i(w)$ will depend only on the left coset
$wW^i$, where  $W^i=W^{\om_i}$ is the centralizer
of $\om_i$; see (\ref{stabwi}) for the list
of $W^i$.

Let  $\ga_w\equal\sum_{i=1}^n a_i \al_i^\vee$ for 
$a_i=a_i(w)=a_i(w\hbox{ mod }W^i)$ as follows:
\begin{align}\label{gammaa}
A_n\,:\
&a_i= \#\{j\le i\,:\, w_j>i\}
\hbox{\ for\ } 
1\le i \le n; \\
\label{gammab}
B_n\,:\  &a_1\!=2 \hbox{\ if\ } w_1=-1,\, a_1=1 \hbox{\ if\ } |w_1|>1
\hbox{\ and\ } 0 \hbox{\ otherwise},\\
&a_{i}\!=\#\{j\le i\,:\, w_j>i\}\!+\!\#\{j\le i\,:\, w_j<0\}
\hbox{\ for\ } 1<i<n,\notag\\ 
\label{gammac}
&a_{n}\!=\bigl[\,(\#\{ 1\le i\le n\,:\, w_j<0\}+1)/2\,\bigr];\notag\\
C_n\,:\
&a_i\!= \#\{j\le i\,:\, w_j>i\}\!+\!\#\{j\le i\,:\, w_j<0\}
\hbox{\ for\ } 
1\le i \le n;\\ 
\label{gammad}
D_n\,:\  &a_1\!=2 \hbox{\ if\ } w_1=-1,\, a_1=1 \hbox{\ if\ } |w_1|>1
\hbox{\ and\ } 0 \hbox{\ otherwise},\\
&a_{i}\!=\#\{j\le i\,:\, w_j>i\}\!+\!\#\{j\le i\,:\, w_j<0\},\, 
1<i<n\!-\!1,\notag\\ 
&a_{n-1}\!=\! \bigl[\,\bigl(\,(1 \hbox{\, if \,} w_j=n 
\hbox{\, for\, } j<n)
+(1 \hbox{\, if\, } 0<w_n<n)+\notag\\
&\ \ \ \ \ \ \ \ \ \ \ (1 \hbox{\, if\, } w_n\!=\!-n)+
\#\{j<n\,:\, -n<w_j<0\}\,\bigr)/2\,\bigr],\notag\\ 
&a_{n}\!=\bigl[\,(\#\{ 1\le i\le n\,:\, w_j<0\})/2\,\bigr].
\notag
\end{align}
\medskip

In the twisted setting (we mark it by $\nu$), let
\begin{align}\label{gammabt}
B^\nu_n\,:\
&a_i = (\,\#\{j\le i\,:\, w_j>i \hbox{ or } w_j<0\}\,)\,\nu_i\,
\hbox{\ for\ } 1\le i \le n; \\
\label{gammact}
C^\nu_n\,:\  &a_1=2 \hbox{ if } w_1=-1,\ a_1=1 \hbox{\ if\ } |w_1|>1
\hbox{\ and\ } 0 \hbox{\ otherwise},\\
&a_{i}\ =\ \#\{j\le i\ :\ w_j>i\, \hbox{\, or\, } \, w_j<0\}
\hbox{\ \ for\ \ } 1<i<n,\notag\\ 
&a_n\ =\, \#\{ 1\le j\le n\,:\, w_j<0\} \hbox{ plus \,} 1 
\hbox {\, if this } \# \hbox{ equals } 1. \notag
\end{align}

\begin{proposition}\label{BCD-GAMMA}
Defining $\ga$  via formulas (\ref{gammaa})- (\ref{gammact}),
one has $\, n(\om_i,w)=a_i(w)=(\ga_{w},\om_i)\, $ 
for $1\le i\le n$, where these numbers depend only
on the left coset $wW^i$.\sq
\end{proposition}

\begin{corollary}\label{ABC-GAMMA}
For the classical root systems
satisfying the total additivity of $n(\la,w)$ (twisted or
untwisted) listed in Theorem \ref{ABDaddit},
one has $\,n(\la,w)=(\ga_{w},\la)$ for any $\la\in P_+$. 
\end{corollary}
\smallskip

\setcounter{equation}{0}
\section{\sc The PBW filtration}
The Kostant $q$-partition function obviously gives the graded 
PBW-characters of the Verma modules (calculated from the 
highest vectors). It is linked to the nil-DAHA and
dag-polynomials as well, but it will be not discussed in 
this paper. We switch in this section to the finite-dimensional 
representations of simple Lie algebras, which is the key in
the {\em PBW-$E^\dag$ correspondence}.

\subsection{\bf General setup}
Given a root system $R$, let $\fg$ be 
the corresponding simple Lie algebra with the Cartan decomposition 
$\fg=\fn\oplus\fh\oplus\fn^-$ and the Weyl group $W$. 
We fix the Borel subalgebra
$\fb=\fh\oplus\fn$ and the Cartan basis
$f_\al$ of $\fn^-$, $\al$ running through the
set of positive roots. For a dominant integral weight
$\la$, let $V_\la$ be the corresponding irreducible $\fg$-module
of highest weight $\la$ with highest weight vector $v_\la$. 
In particular, $V_\la=\rm{U}(\fn^-)v_\la$. The PBW filtration on the
universal enveloping algebra $\rm{U}(\fn^-)$ induces the 
increasing {\em PBW filtration} on $V_\la$:
\[
F_s={\rm span}\{x_1\dots x_lv_\la:\ x_i\in\fn^-, l\le s\}.
\]
The {\em associated graded space} is denoted by $V_\la^a$.
$V_\la^a$ is naturally a cyclic representation of the 
symmetric algebra $S(\fn^-)$ coming from the action of
$\fn^-$, since $f_\al F_s\subset F_{s+1}$, as well a 
representation of the Borel subalgebra $\fb$ because
$\fb$ preserves each $F_s$. These two actions are
combined in the action of the {\em degenerate Lie algebra} 
$\fg^a$ (see \cite{F1}). The spaces $V_\la^a$ are naturally 
graded:
\[
V_\la^a=\bigoplus_{s\ge 0} V_\la^a(s)=F_s/F_{s-1}.  
\]
For a vector $v\in V_\la^a(s)$, we say that the {\em PBW degree} 
of $v$ is equal to $s$; the PBW degree of a vector 
$w\in V_\la$ equals $s$ if $w\in F_s\setminus F_{s-1}$. 
\smallskip

Let us consider $\la=\theta$ as an example. 
The highest weight representation $V_\theta$ 
is isomorphic to
the adjoint representation. In particular, the highest weight
vector $\,v_\theta\,$ is $\,e_\theta\,$ and the lowest
weight vector is $\,f_\theta$. We note that
$ad(f_\theta)^2e_\theta$ is proportional to $f_\theta$. Since
$V_\theta={\rm U}(\fb)f_\theta$ and the PBW filtration
is $\fb$-invariant, the maximal PBW degree of a vector in
$V_\theta$ is two. This agrees with the considerations 
of Section \ref{maxroots}.
\medskip

The following facts about the representations 
$V_\la$ will be used below:
\begin{enumerate}
\item {For any dominant weights $\la$ and $\mu$, there exists an 
embedding of $\fg$-modules 
$V_{\la+\mu}\to V_\la\otimes V_\mu$ sending a highest weight vector 
$v_{\la+\mu}\in V_{\la+\mu}$ to the
tensor product of highest weight vectors $v_\la\otimes v_\mu$.}  
\item {For any $w\in W$, there is only one way to decompose the 
extremal weight $w(\la+\mu)$ into
a sum of a weight of $V_\la$ and a weight of $V_\mu$. Namely, 
this decomposition is nothing but
$w(\la+\mu)=w(\la)+w(\mu)$.} 
\item {The weight subspace of $V_\la$
is one-dimensional for the weight $w(\la)$ and any $w$. 
We fix a vector $v_{w(\la)}$ in this subspace.}   
\end{enumerate}
\smallskip

\subsection{\bf Extremal PBW degree}
We will restrict ourselves to the extremal vectors only, 
which correspond to considering the extremal $E$-dag 
polynomials above. The {\em extremal PBW degree}  
$d(\la,w)$ is the PBW degree of the vector $v_{w(\la)}$
defined in the previous subsection. 

Aiming at the total additivity of $d(\la,w)$, let us
begin with the following inequality.
\begin{lemma}
$d(\la+\mu,w)\ge d(\la,w)+d(\mu,w)$.
\end{lemma}
{\it Proof.}
Let $\gamma_1,\dots,\gamma_N$ be a sequence of 
roots from $R_+$ such that
\[
f_{\gamma_1}\dots f_{\gamma_s}v_{\la+\mu} = v_{w(\la+\mu)}
\hbox{\, in\, } V_{\la+\mu}.
\]

Under the embedding $V_{\la+\mu}\subset V_\la\otimes V_\mu$, 
each $f_\gamma$ is represented as
$f_\gamma\otimes 1+1\otimes f_\gamma$. Since
$v_{w(\la+\mu)}$ is 
represented by $v_{w(\la)}\otimes v_{w(\mu)}$, we obtain
that $N\ge d(w,\la)+d(w,\mu)$.    

\vskip -0.5cm\sq
\smallskip

Let us prove that $d(\la+\mu,w)=d(\la,w)+d(\mu,w)$. 
We will use the following notion of {\em essential collections} due to 
Vinberg; see \cite{V}, \cite{Gor},\cite{F2}. 
First, we order the set of positive 
roots in a sequence $\beta_1,\beta_2,\dots,\beta_N$ in such
a way that if $\beta_i>\beta_j$ then $i<j$. A {\em collection} 
will be a sequence 
$\sigma=(\la;p_{\beta_1},\dots,p_{\beta_N})$, 
where $\la$ is a dominant weight and $p_{\beta_i}\in\Z_{\ge 0}$
(we will call them exponents). 
To such $\sigma$, we attach a vector 
\begin{align}\label{vsigmaf}
v(\sigma)=f_{\beta_1}^{p_1}\dots f_{\beta_N}^{p_N}v_\la\in V_\la.
\end{align}

Second, we introduce a total ordering in the set of collections with 
fixed $\la$ (we only compare collections 
with coinciding $\la$). For a collection $\sigma$, let 
\[
a_k(\sigma)=\sum_{i=k}^N p_i.
\]
For example, $a_1(\sigma)$ is the sum of all 
$p_i$ in (\ref{vsigmaf}). 
Then we order collections lexicographically via $a_i$'s, i.e.
$\sigma>\tau$ if $a_1(\sigma)>a_1(\tau)$ or if
$a_2(\sigma)>a_2(\tau)$ when $a_1(\sigma)=a_1(\tau)$ and  
so on.

\begin{definition}
A collection $(\la,p_{\beta_1},\dots,p_{\beta_N})$ is called
essential if 
$$v(\sigma)\notin \mathrm{span}\{v(\tau):\ \tau<\sigma\}.$$
\end{definition} 
In particular, if $\sigma$ is essential, then the PBW-degree of the 
vector $v(\sigma)$ equals exactly $a_1(\sigma)$, i.e.
is the sum of all exponents $p_{\beta_i}$.

\begin{theorem}\label{pbwadd}
(i) For the component-wise addition of collections
(including $\la$), the essential collections form a semigroup, 
i.e. if  $\sigma$ and $\tau$ are essential then so is $\sigma+\tau$.

(ii) The total additivity for dominant
$\la,\mu\,$ holds:\, $d(\la+\mu,w)= d(\la,w)+d(\mu,w)$.
\end{theorem}
{\it Proof.}
Part $(i)$ is the key here; this is due to Vinberg (see \cite{V}, 
\cite{F2}).
To justify $(ii)$,
let $\sigma$ and $\tau$ be essential collections such that
\[
v(\sigma)=v_{w(\la)}\in V_\la,\ \ v(\tau)=v_{w(\mu)}\in V_\mu.
\] 
Then $\sigma+\tau$ is essential and therefore $v(\sigma+\tau)$ does 
not vanish in the PBW-graded module. This results in
$d(\la+\mu,w)= d(\la,w)+d(\mu,w)$.

\vskip -0.4cm\sq

\subsection{\bf Fundamental modules}
The numbers $d(\la,w)$ are completely determined by the values 
of the PBW degrees in fundamental representations due 
to Theorem \ref{pbwadd}. We will compute such values for 
types $A,C$ (in this subsection), and then for $D,B$ and $G_2$.
\smallskip

{\em Type $A$}. Let $\fg=\msl_n$, $w\in \S_n$ and 
let $\la=\omega_k$ be a fundamental weight. We claim that 
\begin{equation}\label{typeA}
d(\omega_k,w)=\#\{i\le k:\ w(i)>k\}.
\end{equation}
In particular, it gives that $d(\omega_k,w)=n(\omega_k,w)$, 
where $n(\la,w)$ is defined via the $q$-Kostant function. 
Moreover, $d(\la,w)=n(\la,w)$ for all $\la$ and $w$ since
$d(\la,w)$ and $n(\la,w)$ are both additive in $\la$.

Recall that the fundamental module $V_{\omega_k}$ is isomorphic 
to the wedge power 
$\Lambda^k(V)$, where $V$ is the $n$-dimensional vector 
representation of $\msl_n$.
Let $e_1,\dots,e_n$ be the standard basis of $V$. Then the space 
$V_{\omega_k}$ has a basis $e_{\!J}$ labeled by the subsets 
$J\subset \{1,\dots,n\}$, consisting of $k$ elements.
Namely,
\[
e_{\!J}=e_{j_1}\wedge\dots\wedge e_{j_k},\ J=\{j_1,\dots,j_k\}.
\] 
We set 
\[
\deg_k J=\#\{j\in J:\ j>k\}.
\]
It is easy to see that $e_{\!J}\in F_{\deg_k}$, but 
$e_{\!J}\notin F_{\deg_k-1}$   
(i.e. the PBW degree of the vector
$e_{\!J}$ in $V_{\omega_k}$ is exactly $\deg_k J$). 
Now it suffices to use that the extremal vector $v_{w(\omega_k)}$ is
proportional to $e_{w(1),\dots,w(k)}$, which gives \eqref{typeA}.

\begin{corollary}\label{sum}
Let $\la=\sum_{i=1}^{n-1} m_i\omega_i$. Then 
\[
d(\la,w)=\sum_{k:\ w(k)>k} (m_k+\dots + m_{w(k)-1}).
\]
\end{corollary} 
\smallskip

{\em Type $C$}. Let $\fg=\msp_{2n}$, $w\in W$ and 
let $\la=\omega_k$ be a fundamental weight. 
Recall that the Weyl group of type $C$ contains the
permutation group $\S_n$ as well as all sign changes 
$\veps_i\to -\veps_i$. We claim that 
\[
d(\omega_k,w) = \#\{i\le k:\ w(i)>k\} + \#\{i\le k:\ w(i)<0\}.
\]  
The proof can be either deduced from \cite{FFL2} or directly via 
the embedding $\msp_{2n}\subset \msl_{2n}$.
As in type $A$, we obtain that 
$d(\la,w)=n(\la,w)$ for all $\la$ and $w$.
\smallskip

\subsection{\bf Types \texorpdfstring{$D$ and $B$}{D and B}}
For the type $D$,
let $\omega_1,\dots,\omega_n$ be the set of fundamental 
weights of $\mso_{2n}$.
We fix a basis $e_1,\dots,e_{2n}$  of the vector representation 
$V_{\omega_1}$ of $\mso_{2n}$.
In the following we always assume that the orthogonal algebra 
$\mso_{2n}$ 
is defined as the Lie algebra of the Lie  
group leaving invariant the symmetric form on $\C^{2n}$ 
defined by the $2n\times 2n$-matrix in the basis $e_i$:
$$
\left(\begin{array}{ccccc}0 & 0 & 0 & 0 & 1 \\
0 & 0 & 0 & 1 & 0 \\0 & 0 & .\cdot\,{}^\cdot  
& 0 & 0 \\0 & 1 & 0 & 0 & 0 \\ 1 & 0 & 0 & 0 & 0
\end{array}\right).
$$
For a $n\times n$ matrix $A$, let $A^{\tau}$ be the transpose of 
a matrix with respect to the diagonal given by
${i+j}=2n+1$, i.e. 
$A^{\tau}=( a^{\tau}_{i,j})$ with the entries 
$a^{\tau}_{i,j}=a_{2n+1-j,2n+1-i}$ for $A=(a_{i,j}).$
The Lie algebra $\mso_{2n}$ can be then described as the
following set of matrices:
$$
\mso_{2n}=\left\{\left(\begin{array}{cc}
A & B  \\ C & -A^{\tau} \end{array}\right)\big\vert
A,B,C\in Mat_n, \,B=-B^{\tau},\, C=-C^{\tau};
\right\}
$$
with the Cartan subalgebra being 
$\fh=\rm{diag}(t_1,\ldots,t_n,-t_n\ldots,-t_1)$ 
and the Borel subalgebra the upper
triangular matrices in the presentation above.

Recall that $V_{\omega_1}$ is the $2n$-dimensional vector 
representation 
of $\mso_{2n}$ and one has
$V_{\omega_k}\simeq \Lambda^k(V_{\omega_1})$
for $1\le k\le n-2$.
Also, the extremal vectors in $V_{\omega_k}$, $k=1,\dots,n-2$
are the wedge products of the basis vectors $e_i$. 
We have the following proposition:

\begin{proposition}
Let $k=1,\dots,n-2$. Then the PBW degree of the extremal vector 
$v_{w(\omega_k)}$ $\,(w\in W)\,$ equals the $q$-Kostant degree 
$n(\omega_k,w)$ unless there exists a subset $I\subset\{1,\dots,k\}$
such that $\#\{I\}=2s+3$ for $s\ge 0$ and 
\begin{equation}\label{special}
v_{w(\omega_k)}=\bigwedge_{i\le k, i\notin I} e_i\wedge
\bigwedge_{i\in I} e_{2n+1-i}.  
\end{equation}
In the latter case $d(\omega_k,w)=n(\omega_k,w)+1$.
\end{proposition}
{\it Proof.} Using the explicit realization of the orthogonal algebra
given above, one easily checks that generally the 
shortest possible monomial $f_{\beta_1}\dots f_{\beta_m}$ 
such that $\sum\beta_i=\omega_k-w(\omega_k)$\, ($\beta_i$
are positive roots of $\mso_{2n}$) acts nontrivially
on the highest weight vector. Let us show that  
\eqref{special} describes exactly the cases where
we need to use one additional root vector.

Let $I\subset\{1,\dots,k\}$ and      
\[
u=\bigwedge_{i\le k, i\notin I} e_i\wedge\bigwedge_{i\in I} 
e_{2n+i-i}.
\]
To simplify the notation, let $I=\{1,\dots,m\}$
(rename the indices in the general case). Then the 
shortest monomial in terms of $f_\beta$'s changing
the highest weight $\om_k$ (that of $v_{\om_k}$)
to the weight of $\,u\,$ can be represented as follows.
In terms of the standard matrices 
$E_{pq}=(a_{ij}=\delta_{ip}\delta_{jq})$, it is
\begin{align}\label{zerox}
(E_{2n+1-\sigma(1),1}-E_{2n,\sigma(1)})
&(E_{2n+1-\sigma(2),2}-E_{2n-1,
\sigma(2)})\\
&\dots (E_{2n+1-\sigma(m),m}-E_{2n+1-m,\sigma(m)})\notag
\end{align}
for a proper permutation $\sigma\in \S_m$. It is easy to see that  
the result of application of \eqref{zerox} to 
the highest weight vector $e_1\wedge\dots\wedge e_k$ vanishes 
for odd $m$. However, the vector $u$ can be reached in this case 
by using one additional positive root. Namely, 
for $I=\{1,\dots,m\}$ (as above), we begin with applying 
$E_{k+1,1}-E_{2n,2n-k}$ to the 
highest weight vector and then continue using the shortest 
possible sequence (in the $q$-Kostant sense) as in
(\ref{zerox}).

\vskip -0.5cm
\sq

The representations $V_{\omega_{n-1}}$ and $V_{\omega_n}$ are 
exceptional; they are 
{\em spin representations} (see \cite{FH}, Lecture 20). The following 
proposition follows directly from their explicit realization.

\begin{proposition}
The PBW degree of the extremal vectors 
in the spin representations coincide with the q-Kostant degree.
\end{proposition}
\medskip
  
{\em Type $B$}.
The odd orthogonal case $\fg=\mso_{2n+1}$ is parallel to the even one. 
In particular, the fundamental
representations $V_{\omega_k}$, $k=1,\dots,n-1$ are the wedge 
powers of the
vector representation and $V_{\omega_n}$ is the spin representation. 

\begin{proposition}
The PBW degree of an extremal vectors $v_{w(\omega_k)}$ in the 
fundamental
representation $V_{\omega_k}$ of $\mso_{2n+1}$ coincide with the 
$q$-Kostant degree $n(\omega_k,w)$ unless $k=3,\dots,n-1$ and there 
exists a subset 
$I\subset\{1,\dots,k\}$ such that $\#\{I\}=2s+3$ for $s\ge 0$ and 
\[
v_{w(\omega_k)}=\bigwedge_{i\le k, i\notin I} 
e_i\wedge\bigwedge_{i\in I} e_{2n+i-i}.  
\]
In the latter case $d(\omega_k,w)=n(\omega_k,w)+1$.
\end{proposition}
\medskip

\subsection{\bf The system \texorpdfstring{$G_2$}{G2}}
Let $\al_1,\al_2$ be the standard simple roots. The six positive roots
are as follows:
\begin{gather*}
\beta_1=3\al_1+2\al_2,\ \beta_2=3\al_1+\al_2,\ 
\beta_3=2\al_1+\al_2,\\
\beta_4=\al_1+\al_2,\ \beta_5=\al_2,\ \beta_6=\al_1.
\end{gather*}
We note that the ordering of the roots is fixed in such a way that 
if $\beta_i>\beta_j$ then $i<j$. Such choice of ordering is important
for the construction below. 

Let $\la=k\omega_1+l\omega_2$, $k,l\ge 0$. Define the set
$S(\la)\subset\Z_{\ge 0}^6$ consisting of collections
$(s_i)_{i=1}^6$ subject to the relations:
\begin{gather}\label{s5s6}
s_5\le l,\ s_6\le k,\\
s_2+s_3 + s_6\le k+l,\ s_3+s_4+s_6\le k+l,\
s_4+s_5+s_6\le k+l,\notag\\
s_1+s_2+s_3+s_4+s_5\le k+2l,\ s_2+s_3+s_4+s_5+s_6\le k+2l.\notag
\end{gather}
It is proved in \cite{Gor} that the set 
$\prod_{i=1}^6 f_{\beta_i}^{s_i} v_\la$,
$(s_i)_{i=1}^6\in S(\la)$
is a basis of the {\em associate graded space} $V_\la^a$ (see above).
For the fundamental weights we have the following.

\noindent
If $k=1,l=0$, then relations (\ref{s5s6}) reduce to
\[
s_5=0,\ s_6\le 1,\ s_1+s_2+s_3+s_4\le 1,\ s_2+s_3+s_4+s_6\le 1.
\]
If $k=0$, $l=1$, then we have
\begin{gather*}
s_6=0,\ s_2+s_3 \le 1,\ s_3+s_4\le 1,\ s_4+s_5\le 1,\\
s_1+s_2+s_3+s_4+s_5\le 2.
\end{gather*}
Now it is easy to check that the PBW degrees of the extremal vectors 
in $V_{\om_1}^a$ and $V_{\om_2}^a$ coincide with the $q$-Kostant 
degrees from the table in subsection \ref{G2} in the untwisted 
case. 

For instance, let $w=(212121)$. Then 
the difference $\omega_1-w(\omega_1)=4\al_1+2\al_2$ equals
$\beta_1+\beta_6$ and $f_{\beta_1}f_{\beta_6}$ is indeed
an element of our basis for $k=1$, $l=0$. This justifies 
the value $a_1=2$ in the table Section \ref{G2}. 
For $\omega_2$, $\omega_2-w(\omega_2)=6\al_1+4\al_2=2\beta_1$ 
and $f_{\beta_1}^2$ is in the basis for $k=0$, $l=1$, 
which matches $a_2=2$ in this table.
\smallskip

Let us consider the counterexample for the additivity of $q$-Kostant 
extremal degrees for $G_2$ from Section \ref{Cta}.
Let $\,\la=2\omega_1+\omega_2\,$ and $\,w=w_0$. Then
$\,\la-w(\la)=14\al_1+8\al_2=4\beta_1+2\beta_6$. Taking
the values $s_1=4$, $s_6=2$ and $s_2=s_3=s_4=s_5=0$ in 
(\ref{s5s6}) for $k=2$, $l=1$, 
the PBW degree of the extremal vector corresponding to $w_0$ 
equals $6$, while the  $q$-Kostant degree is $5$.
\smallskip

\setcounter{equation}{0}
\section{\sc Extremal dag-degrees}
Let us list the modifications of formulas 
(\ref{gammab}), (\ref{gammad}), (\ref{gammact})
necessary for the extremal dag-degree vs.
$a_i=a_i(w)$ for the Kostant $q$-partition 
function.

\subsection{\bf Classical root systems}
First of all, no modifications
are necessary for the classical root systems
covered by Theorem \ref{ABDaddit}. In the cases of 
untwisted $B_3$ and $D_4$, no modifications are needed 
as well at level of fundamental weights, 
though the numbers $n(\la,w)$  do not satisfy the 
total additivity for these root systems. We set 

\begin{align}\label{gammabdag}
B_n& (n\ge 4)\,:\
\tilde{a}_i = a_i+1 \hbox{\,  for \,} 3\le i\le n-1 
\hbox{\, provided }\\
|w_j|\le i& \hbox{ for } j\le i \hbox{ and }
\#\{j\le i\,:\, \hbox{sgn}(w_j)=-1\}=3+2s (s\in \Z_+);\notag\\
\label{gammacdag}
C^\nu_n& (n\ge 3)\,:\
\tilde{a}^\nu_i\ =\ a_i+1 \hbox{\,\ for \ } 3\le i\le n 
\hbox{\,\,\, \ provided }\\
|w_j|\le i& \hbox{ for } j\le i \hbox{ and }
\#\{j\le i\,:\, \hbox{sgn}(w_j)=-1\}=3+2s (s\in \Z_+);\notag\\
\label{gammaddag}
D_n& (n\ge 5)\,:\
\tilde{a}_i = a_i+1 \hbox{\, for \,} 3\le i\le n-2
\hbox{\, provided }\\
|w_j|\le i& \hbox{ for } j\le i \hbox{ and }
\#\{j\le i\,:\, \hbox{sgn}(w_j)=-1\}=3+2s (s\in \Z_+).\notag
\end{align}
The numbers $\tilde{a}_i(w)$ as well as $a_i(w)$ depend
only on $\,w\hbox{ mod }W^i$, i.e. on the coset $\,wW^i$.
Note that the only difference between these three cases is in
the range $3\le i\le n-1,n,n-2$ of the indices $i$;
the increase from $a_i$ to $\tilde{a}_i$ is always by $1$ in 
the exceptional cases listed above. For the sake of uniformity, 
we will use the notation
$\tilde{a}_i$ for the coefficients $a_i$ that are not included
in (\ref{gammabdag}-\ref{gammaddag}), i.e. remain unchanged.
Sometimes (not always) we will use the notation 
$\tilde{a}^\nu_i$ in the twisted case, as in
(\ref{gammacdag}).
\smallskip

\begin{theorem}\label{DAG-ADDIT}
(i) For any classical root systems under twisted or
untwisted setting, the dag-degrees of for any
$\la\in P_+$ are 
\begin{align*}
&e^\nu(\la,w) = (\tga^\nu_w,\la),\ 
e(\la,w) = (\tga_w,\la), \hbox{\ \, where\ \, }
\tga^{\bep}_w = \sum_{i=1}^n \tilde{a}^{\bep}_i \al^\vee_i
\end{align*} 
for the coefficients $\tilde{a}^{\bep}_i$ 
defined above ($\bep\,$ means $\nu$ or $\varnothing$ in the
twisted or untwisted cases).

(ii) For any classical untwisted root systems, $e(\la,w)=d(\la,w)$
for all $\la\in P_+$ and $w\in W$, which follows from Theorem
\ref{pbwadd} and the calculations of the PBW-degrees for the
fundamental representations performed in the previous sections.
\sq
\end{theorem}
\medskip

\subsection{\bf The case of \texorpdfstring{$G_2$}{G2}}\label{G2}
Next, let us provide  $\ga_w,\ga^\nu_w$ and $\tga_w,\tga^\nu_w$
for $G_2$ in the untwisted and twisted cases.
In the untwisted case, $\ga_w=\tga_w$ for all $w$.
Note that $\ga^\nu_w\neq \tga^\nu_w$ in spite of
the total additivity for $n^\nu(\la,w)$ for twisted $G_2$.
   
The untwisted coefficients $a_i(w)\,(i=1,2)$ 
are given in the $4$th column in the table below;
the twisted coefficients $a^{\nu}_i(w)$ are in the $5$th
(for $\tilde{R}^\nu$) and their tilde-corrections 
$\tilde{a}_i^\nu(w)$ are in the last column.
\smallskip

The elements $w$ from $W$ (the dihedral group of order $12$)
will be presented simply using their reduced decompositions
$w=s_{i_l}\cdots s_{i_1}$, where $l=l(w)$.

We mark the changed values from  $a_i^\nu(w)$ to
$\tilde{a}_i^\nu(w)$ by star; 
such changes affect only $\om_2$ (the
second value) and do not occur for the untwisted $G_2$.
\medskip

The table of $a$-coefficients for $G_2$:
\smallskip

\begin{tabular}{|c||c|c||c|c|c|}  
\hline
$w=s_{i_l}\cdots s_{i_1}$ & $\om_1-w(\om_1)$ & $\om_2-w(\om_2)$
& $a_1a_2$ & $a^\nu_1a^\nu_2$ & $\tilde{a}^\nu_1\tilde{a}^\nu_2$  \\
\hline\hline
id       & $0$ & $0$                              
&  $00$ &  $00$ &  $00\  $\\ \hline
$2     $ & $0$& $\al_2$
&  $01$ &  $03$ &  $03\  $\\ \hline
$1     $ & $\al_1$ & $0$                          
&  $10$ &  $10$ &  $10\  $\\ \hline
$21    $ & $\al_1+\al_2$ & $\al_2$ 
&  $11$ &  $13$ &  $13\  $\\ \hline
$12    $ & $\al_1$ & $3\al_1+\al_2$
&  $11$ &  $12$ &  $13^* $\\ \hline
$212   $ & $\al_1+\al_2$ & $3\al_1+3\al_2$
&  $12$ &  $13$ &  $13\  $\\ \hline
$121   $ & $3\al_1+\al_2$ & $3\al_1+\al_2$
&  $11$ &  $22$ &  $23^* $\\ \hline
$2121  $ & $3\al_1+2\al_2$ & $3\al_1+3\al_2$
&  $12$ &  $23$ &  $23\  $\\ \hline
$1212  $ & $3\al_1+\al_2$ & $6\al_1+3\al_2$
&  $12$ &  $23$ &  $23\  $\\ \hline
$21212 $ & $3\al_1+2\al_2$ & $6\al_1+4\al_2$
&  $12$ &  $24$ &  $26^* $\\ \hline
$12121 $ & $4\al_1+2\al_2$ & $6\al_1+3\al_2$
&  $22$ &  $23$ &  $23\  $\\ \hline
$212121$ & $4\al_1+2\al_2$ & $6\al_1+4\al_2$
&  $22$ &  $24$ &  $26^* $\\ \hline
\end{tabular}
\medskip

We arrive at the following proposition.

\begin{proposition}\label{GAG2}
Using the table above in the case of $G_2$, for
\begin{align*}
&\ga_w=\tga_w=\sum_{i=1}^2 a_i(w)\al_i^\vee,\ \, 
\ga^\nu_w=\sum_{i=1}^2 a^\nu_i(w)\al_i^\vee,\ \,
\tga^\nu_w=\sum_{i=1}^2 \tilde{a}^\nu_i(w)\al_i^\vee,
\end{align*}
\begin{align}\label{nomg2}
&n(\om_i,w)=(\ga_w,\om_i) \hbox{\, and\, } 
n^\nu(\om_i,w)=(\ga^\nu_w,\om_i)\hbox{\, for\, } i=1,2,\\
\label{nbeg2} 
&e(\la,w)\ =(\tilde{\ga}_w,\la) \hbox{\, and\  } 
e^\nu(\la,w)=(\tilde{\ga}^\nu_w,\la)\hbox{\, for any \, } \la\in P_+. 
\end{align}
For untwisted $G_2$, 
$e(\la,w)= d(\la,w)$ for any $\la\in P_+,w\in W$.
\sq
\end{proposition}
\medskip

\comment{
{\bf Comment}.
Let us give a possible {\em twisted PBW-interpretation} of
the increases marked by $\ast$ in the table above. 
The {\em twisted}
$\mathfrak{g}^\nu$ will be defined as 
$\Bigl(\oplus_{\al\in R}\,\mathfrak{g}_{\al}\,\Bigr)
\oplus \mathfrak{h}$
with the commutators as from $\mathfrak{g}$ of type
$G_2$ subject to the following change: 
\begin{align}\label{xalxbe}
[x_{\al},x_{\be}]\,=\,0 \hbox{\, for short\, } \al\in R,
\hbox{\, where\, } x_\al\in g_\al.
\end{align}

Let us assume that the {\em twisted PBW-$E^\dag$ correspondence},
holds for $\mathfrak{g}^\nu$.
Here each $f_\al$ for long $\al\in R_+$ 
is counted as $3\,(=\nu_{\lng})$. 
The latter alone is insufficient; $\mathfrak{g}$ 
must be changed to $\mathfrak{g}^\nu$ as well. Then the corrections 
marked by ${}^\ast$ in the table lead to the following conclusions.

Let $V_{\om_2}$ be the irreducible representation
of $\mathfrak{g}^\nu$ with the highest vector
$v$ of weight $\om_2$ with respect to 
$\mathfrak{b}^\nu_+\supset \mathfrak{n}_+^\nu=
\oplus_{\al>0} \C e_{\al}.$
Then the increase from 
$24$ to $26^*$, corresponding to taking 
$w_0=212121=s_2\cdots s_2 s_1$, 
means that
\begin{align}\label{g2-26}
&f_{\al_1+\al_2}f_{2\al_1+\al_2}f_{3\al_1+2\al_2}v=0=
f_{\al_1+\al_2}f_{2\al_1+\al_2}f_{\al_1+\al_2}f_{2\al_1+\al_2}v.
\end{align}
Here changing the order of $f$-operators does not matter since
they are all pairwise commutative (in $\mathfrak{g}^\nu$).
Each of these elements must be proportional to the lowest vector 
in $V_{\om_2}$. However they have the twisted PBW degree $5$ or 
smaller and therefore must vanish if the twisted
{\em PBW-$E^\dag$} correspondence holds.
 

Similarly, the increase from $22$ to $23^*$ (equivalently,
from $12$ to $13^*$) means that 
\begin{align}\label{g2-13}
&f_{\al_1}f_{2\al_1+\al_2} v=0=f_{2\al_1+\al_2}f_{\al_1} v.
\end{align}

We will not discuss here a systematic representation theory
of twisted $\mathfrak{g}$ (which is actually part of the
theory of twisted Kac-Moody algebras). In this particular
case, $\om_2=\th=3\al_1+2\al_2$ and the $\mathfrak{g}^\nu$-module
$V_{\om_2}$ is the {\em twisted 
adjoint representation} with $v=e_{\th}$. Thus
the second relation in (\ref{g2-13}) is obvious in the
twisted or untwisted case because $\,(\al_1,\th)=0$. 
The first holds only in the twisted setting; use that
$\,\al_1+\al_2=\th-(2\al_1+\al_2)\,$ is short. Finally,
the PBW-degree is $3$ here due to $\,f_{3\al_1+\al_2}v\neq 0$.

Similarly, (\ref{g2-26}) holds
due to the fact that $\,\al_1+\al_2,\, 2\al_1+\al_2\,$ are
both short, and $\,f_\th^2 v\neq 0\,$ gives that the
twisted PBW-degree is $6$ in this case.

The remaining two cases are $\om_2-w(\om_2)=3\al_1+3\al_2$
and $\om_2-w(\om_2)=6\al_1+3\al_2$. 
The theory of the twisted PBW-filtration will be a subject
of our further research.
\sq
}

Let us provide (computer-generated)
formulas for the deviations of the dag-degrees for $F_4$ and $E_6$  
vs.\! those calculated on the basis of the Kostant
$q$\~partition function. We set $A[i]=X_{\al_i}$ for $1\le i\le n$;
accordingly, 
$$
A[c]=A[c_1,\ldots,c_n]=\prod_{i=1}^n \,A[i]^{\,c_i} \hbox{\, for\, }
Q\ni c=\sum_{i=1}\, c_i \al_i.
$$
For a given fundamental weight
$\om_i$, let 
\begin{align}\label{Pomwom}
\tilde{\E}^\dag_i=\sum_{w\in W/W^i}A[\om_i-w(\om_i)]q^{-e(\om_i,w)},
1\le i\le n, 
\end{align}
where $W^i=W^{\om_i}$ is the centralizer
of $\om_i$. Recall that the dag-degree $e(\la,w)$ depends
only on $w(\la)$, so the summation here is over $W/W^i$.

Only {\em singular  monomials\,} in $\,\tilde{\E}^\dag_i\,$
are considered in the following two subsections i.e. those with 
$q$-degrees different from the corresponding ones obtained via the 
Kostant $q$-partition function. The notation
will be $\,\tilde{\E}^\ddag_i\,$ for such singular subsums\,
($1\le i\le n$).  

The singular monomials are of obvious interest
since they describe (proven for $ABCD$ and $G$) 
the extremal weights $\la$ in the untwisted $\mathfrak{g}$-modules 
$V_{\om_i}\ni v=$vac\, such that 
$\prod_{\al}f_\al(v)=0$ for the shortest possible sequence
$\{\al\in R_+\}$ satisfying $\om_i-\la=\sum \al$. 

Accordingly, we provide the corresponding singular extremal terms of
the Kostant $q$-partition function, namely the 
polynomials
\begin{align}\label{Pomwomka}
\K^{\sing}_i=\!\!\sum_{w\in W/W^i}A(\om_i-w(\om_i))q^{-n(\om_i,w)}
\hbox{ s.t. } e(\om_i,w)>n(\om_i,w).
\end{align}
All monomials with $e(\om_i,w)\neq
n(\om_i,w)$ are counted in $\K^{\sing}_i$, since the inequality
can be only in this direction:\, $e(\om_i,w)>n(\om_i,w)$.
\medskip

\subsection{\bf The system \texorpdfstring{$F_4$}{F4}}

We begin with the {\em twisted\,} singular sums $\tilde{\E}^\ddag\,$ 
for $F_4$; to avoid any misunderstanding we put $\,\nu\,$.
The results of our computer calculations are as follows:
\smallskip

$(\tilde{\E}^\ddag_1)^\nu=$
\renewcommand{\baselinestretch}{0.5} 
{\small
\noindent
\(
\,\frac{A[3,3,4,2]}{q^4}
+\frac{A[3,4,4,2]}{q^4}+\frac{A[3,4,6,2]}{q^4}+\frac{A[3,4,6,4]}{q^4}
+\frac{A[3,5,6,2]}{q^4}+\frac{A[3,5,6,4]}{q^4}
+\frac{A[3,5,8,4]}{q^4}+\frac{A[3,6,8,4]}{q^4},
\)
}
\renewcommand{\baselinestretch}{1.2} 
\smallskip

$(\K^{\sing}_1)^\nu=$
\renewcommand{\baselinestretch}{0.5} 
{\small
\noindent
\(
\frac{A[3,3,4,2]}{q^3}+\frac{A[3,4,4,2]}{q^3}+\frac{A[3,4,6,2]}{q^3}
+\frac{A[3,4,6,4]}{q^3}+\frac{A[3,5,6,2]}{q^3}+\frac{A[3,5,6,4]}{q^3}
+\frac{A[3,5,8,4]}{q^3}+\frac{A[3,6,8,4]}{q^3};\ \ 
\)
}
\renewcommand{\baselinestretch}{1.2} 

\noindent
\smallskip

$(\tilde{\E}^\ddag_2)^\nu=$
\renewcommand{\baselinestretch}{0.5} 
{\small
\noindent
\(
\frac{A[0,3,4,2]}{q^4}+\frac{A[1,3,4,4]}{q^4}+\frac{A[1,3,6,2]}{q^4}
+\frac{A[1,3,6,4]}{q^4}+\frac{A[2,3,4,4]}{q^4}+\frac{A[2,3,6,2]}{q^4}
+\frac{A[2,3,6,4]}{q^4}+\frac{A[2,7,8,4]}{q^6}+\frac{A[2,7,10,4]}{q^6}
+\frac{A[2,7,10,6]}{q^6}+\frac{A[3,3,4,2]}{q^4}+\frac{A[3,7,8,2]}{q^6}
+\frac{A[3,7,8,6]}{q^6}+\frac{A[3,7,12,6]}{q^6}+\frac{A[4,7,8,2]}{q^6}
+\frac{A[4,7,8,6]}{q^6}+\frac{A[4,7,12,6]}{q^6}
+\frac{A[4,9,10,4]}{q^6}
+\frac{A[4,9,10,6]}{q^6}+\frac{A[4,9,12,4]}{q^6}
+\frac{A[4,9,12,8]}{q^6}
+\frac{A[4,9,14,6]}{q^6}+\frac{A[4,9,14,8]}{q^6}
+\frac{A[5,7,8,4]}{q^6}
+\frac{A[5,7,10,4]}{q^6}+\frac{A[5,7,10,6]}{q^6}
+\frac{A[5,9,10,4]}{q^6}
+\frac{A[5,9,10,6]}{q^6}+\frac{A[5,9,12,4]}{q^6}
+\frac{A[5,9,12,8]}{q^6}
+\frac{A[5,9,14,6]}{q^6}+\frac{A[5,9,14,8]}{q^6}
+\frac{A[5,11,14,6]}{q^8}
+\frac{A[5,11,14,8]}{q^8}+\frac{A[5,11,16,8]}{q^8}
+\frac{A[6,11,14,6]}{q^8}+\frac{A[6,11,14,8]}{q^8}
+\frac{A[6,11,16,8]}{q^8}+\frac{A[6,12,16,8]}{q^8}, \ \ 
\)
}
\renewcommand{\baselinestretch}{1.2} 
\smallskip

$(\K^{\sing}_2)^\nu\!=\!$
\renewcommand{\baselinestretch}{0.5} 
{\small
\noindent
\(
\frac{A[0,3,4,2]}{q^3}+\frac{A[1,3,4,4]}{q^3}+\frac{A[1,3,6,2]}{q^3}
+\frac{A[1,3,6,4]}{q^3}+\frac{A[2,3,4,4]}{q^3}+\frac{A[2,3,6,2]}{q^3}
+\frac{A[2,3,6,4]}{q^3}+\frac{A[2,7,8,4]}{q^5}+\frac{A[2,7,10,4]}{q^5}
+\frac{A[2,7,10,6]}{q^5}+\frac{A[3,3,4,2]}{q^3}+\frac{A[3,7,8,2]}{q^5}
+\frac{A[3,7,8,6]}{q^5}+\frac{A[3,7,12,6]}{q^5}+\frac{A[4,7,8,2]}{q^5}
+\frac{A[4,7,8,6]}{q^5}+\frac{A[4,7,12,6]}{q^5}
+\frac{A[4,9,10,4]}{q^5}
+\frac{A[4,9,10,6]}{q^5}+\frac{A[4,9,12,4]}{q^5}
+\frac{A[4,9,12,8]}{q^5}
+\frac{A[4,9,14,6]}{q^5}+\frac{A[4,9,14,8]}{q^5}
+\frac{A[5,7,8,4]}{q^5}
+\frac{A[5,7,10,4]}{q^5}+\frac{A[5,7,10,6]}{q^5}
+\frac{A[5,9,10,4]}{q^5}
+\frac{A[5,9,10,6]}{q^5}+\frac{A[5,9,12,4]}{q^5}
+\frac{A[5,9,12,8]}{q^5}
+\frac{A[5,9,14,6]}{q^5}+\frac{A[5,9,14,8]}{q^5}
+\frac{A[5,11,14,6]}{q^6}
+\frac{A[5,11,14,8]}{q^6}+\frac{A[5,11,16,8]}{q^6}
+\frac{A[6,11,14,6]}{q^6}+\frac{A[6,11,14,8]}{q^6}
+\frac{A[6,11,16,8]}{q^6}+\frac{A[6,12,16,8]}{q^6}; \ \ 
\)
}
\renewcommand{\baselinestretch}{1.2} 

\smallskip

$(\tilde{\E}^\ddag_3)^\nu=$
\renewcommand{\baselinestretch}{0.5} 
\renewcommand{\baselinestretch}{0.5} 
{\small
\noindent
\(
\frac{A[1,2,3,3]}{q^3}+\frac{A[3,7,11,5]}{q^5}+\frac{A[3,7,11,6]}{q^5}
+\frac{A[4,7,11,5]}{q^5}+\frac{A[4,7,11,6]}{q^5}+
\frac{A[4,8,11,5]}{q^5}
+\frac{A[4,8,11,6]}{q^5}+\frac{A[4,8,12,6]}{q^6},\ \ 
\)
}
\renewcommand{\baselinestretch}{1.2}
\smallskip

$(\K^{\sing}_3)^\nu\!=\!$
\renewcommand{\baselinestretch}{0.5} 
{\small
\noindent
\(
\frac{A[1,2,3,3]}{q^2}+\frac{A[3,7,11,5]}{q^4}
+\frac{A[3,7,11,6]}{q^4}
+\frac{A[4,7,11,5]}{q^4}+\frac{A[4,7,11,6]}{q^4}+
\frac{A[4,8,11,5]}{q^4}
+\frac{A[4,8,11,6]}{q^4}+\frac{A[4,8,12,6]}{q^4}.\ \ 
\)
}
\renewcommand{\baselinestretch}{1.2} 

\smallskip

Finally, 
$(\tilde{\E}^\ddag_4)^\nu=0=(\K^{\sing}_4)^\nu.$
\medskip

{\em Untwisted $F_4$.}
\smallskip

One has: \ $\tilde{\E}^\ddag_1=0=\K^{\sing}_1;$ 

$\tilde{\E}^\ddag_2=$
\renewcommand{\baselinestretch}{0.5} 
{\small
\noindent
\(
\frac{A[3,3,4,2]}{q^3}+\frac{A[5,11,14,6]}{q^5}
+\frac{A[5,11,14,8]}{q^5}+\frac{A[5,11,16,8]}{q^5}
+\frac{A[6,11,14,6]}{q^5}+\frac{A[6,11,14,8]}{q^5}
+\frac{A[6,11,16,8]}{q^5}+\frac{A[6,12,16,8]}{q^6},\ \ 
\)
}
\renewcommand{\baselinestretch}{1.2} 
\smallskip

$\K^{\sing}_2=$
\renewcommand{\baselinestretch}{0.5} 
{\small
\noindent
\(
\frac{A[3,3,4,2]}{q^2}+\frac{A[5,11,14,6]}{q^4}
+\frac{A[5,11,14,8]}{q^4}+\frac{A[5,11,16,8]}{q^4}
+\frac{A[6,11,14,6]}{q^4}+\frac{A[6,11,14,8]}{q^4}
+\frac{A[6,11,16,8]}{q^4}+\frac{A[6,12,16,8]}{q^4};\ \ 
\)
}
\renewcommand{\baselinestretch}{1.2} 

\smallskip

$\tilde{\E}^\ddag_3=$
\renewcommand{\baselinestretch}{0.5} 
{\small
\noindent
\(
\frac{A[3,7,11,5]}{q^4}+\frac{A[3,7,11,6]}{q^4}
+\frac{A[4,7,11,5]}{q^4}+\frac{A[4,7,11,6]}{q^4}
+\frac{A[4,8,11,5]}{q^4}
+\frac{A[4,8,11,6]}{q^4}+\frac{A[4,8,12,6]}{q^4},\ \ 
\)
}
\renewcommand{\baselinestretch}{1.2} 
\smallskip

$\K^{\sing}_3=$
\renewcommand{\baselinestretch}{0.5} 
{\small
\noindent
\(
\frac{A[3,7,11,5]}{q^3}+\frac{A[3,7,11,6]}{q^3}
+\frac{A[4,7,11,5]}{q^3}+\frac{A[4,7,11,6]}{q^3}
+\frac{A[4,8,11,5]}{q^3}
+\frac{A[4,8,11,6]}{q^3}+\frac{A[4,8,12,6]}{q^3}.\ \ 
\)
}
\renewcommand{\baselinestretch}{1.2}

\smallskip

Finally, $\tilde{\E}^\ddag_4=0=\K^{\sing}_4.$
See the Appendix for the complete list of
(computer-generated)  $\,\tilde{\E}_{i}^\dag\,$
for the twisted and untwisted $F_4$.
\medskip

\subsection{\bf Singular monomials for \texorpdfstring{$E_6$}{E6}}
We will omit in this paper the complete list of 
$\tilde{\E}^\dag$-polynomials 
in this case (known for $E_{6,7}$) and will 
provide only the subsums of the {\em singular\, }
monomials, $\tilde{\E}^\ddag$,
i.e. those with $q$-degrees different
from the corresponding ones obtained via the 
Kostant $q$-partition function. 

First of all,
$\,\tilde{\E}^\ddag_i\,=\,0\,=\, \K^{\sing}_i\,$
for $\,i=1,2,6$. For the other $\,i\,$, the results
of our calculations are as follows.

\smallskip

$\tilde{\E}^\ddag_3=$
\renewcommand{\baselinestretch}{0.5} 
{\small
\noindent
\(
\frac{A[3,3,6,6,4,2]}{q^4}+\frac{A[3,3,6,7,4,2]}{q^4}+
\frac{A[3,3,6,7,5,2]}{q^4}+\frac{A[3,3,6,7,5,3]}{q^4}+
\frac{A[3,4,6,7,4,2]}{q^4}+\frac{A[3,4,6,7,5,2]}{q^4}+
\frac{A[3,4,6,7,5,3]}{q^4}+\frac{A[3,4,6,8,5,2]}{q^4}+
\frac{A[3,4,6,8,5,3]}{q^4}+\frac{A[3,4,6,8,6,3]}{q^4}.\ \
\)
}
\renewcommand{\baselinestretch}{1.2} 

The corresponding $\K^{\sing}_3$ (with the Kostant $q$-degrees)
results from $\tilde{\E}^\ddag_3$ when all
$q$-degrees in the denominators are diminished from $4$
to $3$.
\smallskip

$\tilde{\E}^\ddag_4=$
\renewcommand{\baselinestretch}{0.5} 
{\small
\noindent
\(
\frac{A[0,3,3,6,4,2]}{q^4}+\frac{A[1,3,2,3,2,1]}{q^3}
+\frac{A[2,3,4,6,3,0]}{q^4}+\frac{A[2,3,4,6,3,3]}{q^4}
+\frac{A[2,3,4,6,6,3]}{q^4}+\frac{A[2,5,6,10,6,3]}{q^5}
+\frac{A[2,5,6,10,7,3]}{q^5}+\frac{A[2,5,6,10,7,4]}{q^5}
+\frac{A[3,3,3,6,4,2]}{q^4}+\frac{A[3,3,6,6,4,2]}{q^4}
+\\ \frac{A[3,5,6,10,6,2]}{q^5}+\frac{A[3,5,6,10,6,4]}{q^5}
+\frac{A[3,5,6,10,8,4]}{q^5}+\frac{A[3,5,7,10,6,2]}{q^5}
+\frac{A[3,5,7,10,6,4]}{q^5}+\\ \frac{A[3,5,7,10,8,4]}{q^5}
+\frac{A[3,5,7,11,7,3]}{q^5}+\frac{A[3,5,7,11,7,4]}{q^5}
+\frac{A[3,5,7,11,8,4]}{q^5}+\frac{A[3,6,7,11,7,3]}{q^5}
+\\ \frac{A[3,6,7,11,7,4]}{q^5}+\frac{A[3,6,7,11,8,4]}{q^5}
+\frac{A[4,5,6,10,6,3]}{q^5}+\frac{A[4,5,6,10,7,3]}{q^5}
+\frac{A[4,5,6,10,7,4]}{q^5}+\\ \frac{A[4,5,7,10,6,2]}{q^5}
+\frac{A[4,5,7,10,6,4]}{q^5}+\frac{A[4,5,7,10,8,4]}{q^5}
+\frac{A[4,5,7,11,7,3]}{q^5}+\frac{A[4,5,7,11,7,4]}{q^5}
+\\ \frac{A[4,5,7,11,8,4]}{q^5}+\frac{A[4,5,8,10,6,3]}{q^5}
+\frac{A[4,5,8,10,7,3]}{q^5}+\frac{A[4,5,8,10,7,4]}{q^5}
+\frac{A[4,5,8,11,7,3]}{q^5}+\\ \frac{A[4,5,8,11,7,4]}{q^5}
+\frac{A[4,5,8,11,8,4]}{q^5}+\frac{A[4,6,7,11,7,3]}{q^5}
+\frac{A[4,6,7,11,7,4]}{q^5}+\frac{A[4,6,7,11,8,4]}{q^5}+\\ 
\frac{A[4,6,8,11,7,3]}{q^5}+\frac{A[4,6,8,11,7,4]}{q^5}
+\frac{A[4,6,8,11,8,4]}{q^5}+\frac{A[4,6,8,12,8,4]}{q^6}.\ \
\)
}
\renewcommand{\baselinestretch}{1.2} 

The corresponding $\K^{\sing}_4\,$ is obtained when
$q^5$ is replaced by $q^4$ in the singular monomials
with the following exceptions.
The singular monomial $A[4,6,8,12,8,4]/q^6\,$ must be replaced by
$A[4,6,8,12,8,4]/q^4\,$ and  $A[1,3,2,3,2,1]/q^3\,$ 
by  $A[1,3,2,3,2,1]/q^2$. For the monomials
\begin{align*}
&A[3,3,6,6,4,2],\ A[3,3,3,6,4,2],\ A[2,3,4,6,6,3],\\ 
&A[2,3,4,6,3,3],\ A[2,3,4,6,3,0],\ A[0,3,3,6,4,2],
\end{align*}
the power $q^4$ must be replaced by $q^3$. 
\smallskip

$\tilde{\E}^\ddag_5=$
\renewcommand{\baselinestretch}{0.5} 
{\small
\noindent
\(
\frac{A[2,3,4,6,6,3]}{q^4}+\frac{A[2,3,4,7,6,3]}{q^4}
+\frac{A[2,3,5,7,6,3]}{q^4}+\frac{A[2,4,4,7,6,3]}{q^4}
+\frac{A[2,4,5,7,6,3]}{q^4}+\frac{A[2,4,5,8,6,3]}{q^4}
+\frac{A[3,3,5,7,6,3]}{q^4}+\frac{A[3,4,5,7,6,3]}{q^4}
+\frac{A[3,4,5,8,6,3]}{q^4}+\frac{A[3,4,6,8,6,3]}{q^4}.
\)
}
\renewcommand{\baselinestretch}{1.2} 

Here $\K^{\sing}_5$ is obtained from $\tilde{\E}^\ddag_5\,$ when
$q^4$ in the denominators is replaced by $q^3$.
\medskip

\setcounter{equation}{0}
\section{\sc Generalizations, perspectives}\label{SEC:GENPER}
Following \cite{CO1}, we will provide in this section
the definition of the full dag-polynomials 
(not only their extremal parts). We begin with 
the general nonsymmetric Macdonald polynomials $E_b(X;q,t)$;
see \cite{Op,Ma,Ch1}. The twisted and untwisted 
settings will be considered; accordingly, we set 
$\bep\,=\nu$ and $\bep\,=\varnothing$ as above. 
We note that paper \cite{CO1} is written in the 
twisted setting, but the untwisted case is parallel.
\smallskip

\subsection{\bf Nonsymmetric polynomials}
The affine {\em Demazure-Lusztig operators} are
\begin{align}
&T_i\  = \  t^{1/2} s_i\ +\
(t^{1/2}-t^{-1/2})(X_{\al_i}-1)^{-1}(s_i-1),
\ 0\le i\le n;
\label{Demazx}
\end{align}
they obviously preserve $\Z[q,t^{\pm 1/2}][X_b,\,b\in P]$.
We note that only the formula for $T_0$ involves $q$\,:
\begin{align}
&T_0\,  = \, t^{1/2}s_0\ +\ (t^{1/2}-t^{-1/2})
(X_0 -1)^{-1}(s_0-1),\hbox{\, where\ }\notag\\
&X_0=qX_{\th^\bep}^{-1},\
s_0(X_b)\ =\ X_bX_{\th^\bep}^{-(b,\th^\bep)}
 q^{(b,\th^\bep)},\
\al_0\ =\ [-\th^\bep,1].
\end{align}

For $\hat{W}^\bep\ni \hw =\pi_r s_{i_l}\cdots s_{i_1}$,
where $l=l(\hw)$, the element  
 $T_{\hw}=\pi_r s_{i_l}\cdots s_{i_1}$ does not depend on
the choice of the reduced decomposition of $\hw$.
For the sake of uniformity, let $\check{P}^\bep$ be $P$
in the twisted case (for $\bep\,=\nu$) and  
$\check{P}^\bep=P^\vee$ in the untwisted case ($\bep\,=\varnothing$).

We set $Y_b=T_b$ for $b\in \check{P}_+^\bep$. Then 
$Y_{b-c}=Y_b Y_c^{-1}$  depends only on $b-c$ and 
(following Bernstein-Zelevinsky-Lusztig) this 
can be used to define $Y_b$ for any $b\in P^\bep$.
We set $\rho^\bep=\rho$ for $\bep\,=\varnothing$ and
$\rho^\bep=\check{\rho}$ for $\bep\,=\nu$.

For generic parameters $q,t$, the nonsymmetric polynomials
(also called $E$-polynomials)  $\{E_b=E_b(X;q,t), b\in P\} 
\subset \Q(q,t^{1/2})[X_b, b\in P]$ can be defined as 
eigenfunctions of the operators
$Y_a$ ($a\in \check{P}^\bep$). This fixes them uniquely
up to proportionality. More explicitly, 
\begin{align}
&Y_a(E_b)\ =\ q^{-(a,b)}\,t^{(a,u_b^{-1}(\rho^{\bep}))}\,E_b 
\hbox{\ \ for\ \ }  
a\in \check{P}^\bep, 
\label{Yone}
\end{align}
where $u_b\in W$ is the element of minimal 
length such that $u_b(b)\in P_-$ (it is unique). Equivalently,
$b=\pi_b u_b$ such that $l(b)=l(\pi_b)+l(u_b)$ and 
$u_b\in W$ is of maximal possible length. This definition
extends  (\ref{xwo}); see, e.g. Proposition 1.2 from \cite{CO1}.

For $b\in P$, let $b_-\equal u_b(b),\ b_+\equal w_0u_b(b).$
We will also use $\,\iota(b)=b^\iota=-w_0(b)\,$ and will extend this 
automorphism to $\R^{n+1}$
as follows:\, $\iota([z,\ze])=-w_0([z,\ze])=[-w_0(z),\ze]$.
In particular, $\iota(\al_i)=\al_{i^\iota}$ for the image
$i^\iota$ of $i=0,1,\ldots, n\,$ under the action of the
automorphism $-w_0$ extended to the
completed Dynkin diagram $\tilde{\Ga}^\bep$
by the relation $\al_0^\iota=\al_0$. We naturally set 
$s_i^\iota=s_{i^\iota}$ and $\pi_r^\iota=\pi_{r^\iota}$ 
for $i\ge 0,r\in O$. For instance, $\pi_b^\iota=\pi_{b^\iota}\,$, 
$u_b^\iota=u_{b^\iota}\,$.

One has:
\begin{align}
&E_b\in \oplus_c\Q(q,t)X_a\hbox{\, where\, }
a_-\in b_-+Q_+ \hbox{\ \,and if\ }\label{macdfilt}\\
&a_-=b_-\hbox{ then }  u_a \geq u_b
\hbox{\, for the Bruhat order\, }\geq. 
\label{macdbruhat}
\end{align}
Let $a\succeq b\,$ if $\,u_a\geq u_b;\,$ we will use
$\succ$ when the inequality is strict in the Bruhat
order. The polynomials $E_b$ are normalized by the 
condition $E_b-X_b\in \oplus_{a\succ b}\Q(q,t)X_a$.
\smallskip

Generalizing formula (2.54) from \cite{CO1},
\begin{align}
E_b^\star &=q^{-(c,b^\iota)}
t^{l(u_b)-l(w_0)/2+(c,u_{b^\iota}^{-1}(\rho^\bep))}
T_{w_0}Y_c^{-1}(E_{b^\iota}),
\label{E-ast-rengen}
\end{align}
for any $c\in P^\bep$,
where $X_b^\star=X_b^{-1}, q^\star=q^{-1},t^\star=t^{-1}$.
Recall that $\rho^\nu=\check{\rho}$ (when $\,\bep\,=\nu$, 
i.e. in the twisted case). 
\medskip

\subsection{\bf Bar- and dag-polynomials}
We define them as
follows:
\begin{align}\label{Ebardag}
\overline{E}_b(X;q)\! =\! E_b(X;q,t\to 0),\ 
E^\dag_b(X;q)\! =\! E_b(X;q,t\to\infty),\, b\in P.
\end{align} 
See \cite{CO1} for the justification of their existence
(minor modifications are needed in the untwisted case).
Here the theory of DAHA and nil-DAHA is the foundation 
(as well as in the general theory of $E$-polynomails).
 
We will use $\overline{T}_i\equal (t^{1/2}T_i)\mid_{t\to 0}$
and $\overline{T}_i^\prime\equal (t^{1/2}T_i^{-1})\mid_{t\to 0}=
\overline{T}_i+1$.
Then (\ref{E-ast-rengen}) results in the
following generalization of (2.50) from
\cite{CO1} to the case of arbitrary $b$ (not only
antidominant). 

\begin{proposition}\label{ASTBARDAG}
For $b\in P$, let $c\in P^\bep$ be any element such
that
\begin{align}\label{c-conds}
u_c\!=\! w_0u_{b^\iota}\!=\! u_{b}w_0,
\hbox{\, i.e. } u_{b^\iota} \hbox{ is maximal }\, u\,
\hbox{ satisfying } u(c)\in P_+. 
\end{align}
Then given a reduced decomposition
$\pi_c=\pi_r s_{j_l}\cdots s_{j_1}\, (r\in O)$, one has
\begin{align}
E_b^\dag\ &=\ q^{-(c_+,b_+)}\,
\bigl(\overline{T}_{u_{b}^{-1}}
\overline{T}_{j_1}^{\,\ep_1}\overline{T}_{j_2}^{\,\ep_2}\cdots
\overline{T}_{j_l}^{\,\ep_l}\,\pi_r^{-1}
\ (\overline{E}_{b^\iota})\bigr)^\ast\notag\\
&=\ q^{-(c_+,b_+)}\,
\bigl(\,\pi_{r^\iota}\ \overline{T}_{j_l^\iota}^{\,\ep_l}\cdots
\,\overline{T}_{j_2^\iota}^{\,\ep_2}
\overline{T}_{j_1^\iota}^{\,\ep_1}\ \overline{T}_{u_{b}^\iota}
\ (\overline{E}_{b^\iota})\bigr)^\ast,
\label{ebminusdag}
\end{align}
where $(q^mX_b)^\ast\!=\!q^{-m}X_{-b}$ and 
$\{\ep_p\}$ are defined as follows in terms of the
sequence 
\begin{align}\label{alsequ}
\al^1=\al_{j_1},\,\al^2=s_{j_1}(\al_{j_2}),\,
\al^3=s_{j_1}s_{j_2}(\al_{j_3}),\hbox{\, and so on,}
\end{align}
from (\ref{alsuperseq}) representing $\la(\pi_c)$.
We pick $\ep_p=\prime$ 
if the nonaffine component of $u_b^{-1}(\al^p)$ is negative 
and $\ep_p=\varnothing$ otherwise.
\end{proposition}
{\it Proof.}
See formula (1.42) from \cite{CO1}  concerning picking 
$\{\ep_p\}$. The relation $u_c\ =\ u_b w_0\,$ gives that
$$
u_c(c)=c_-=w_0 u_{b^\iota}(c)\hbox{\,\, and therefore\,\, }
c=u_{b^\iota}^{-1}w_0(c_-).
$$
Using  $w_0(c_-)=c_+$ and
$b^\iota=u_{b^\iota}^{-1}(-b_+)$, one
obtains that $(c,b^\iota)=-(c_+,b_+)$. Similarly,
$(c,u_{b^\iota}^{-1}(\rho^\bep))=$ 
$(c_+,\rho^\bep).$ Thus  (\ref{E-ast-rengen})
reads as
\begin{align}
E_b^\star &=q^{(c_+,b_+)}
t^{l(u_b)-l(w_0)/2+(c_+,\rho^\bep))}
T_{w_0}Y_c^{-1}(E_{b^\iota}).
\label{E-ast-rengena}
\end{align}
This power of $\,t\,$ must be added to $T_{w_0}Y_c^{-1}$
to ensure the existence of the limit  $t\to 0$. Indeed, 
$t^{l(w_0)+(c^+,\rho^\bep)}$
is needed if there is no reduction in $T_{w_0}Y_c^{-1}$;
however, $T_{w_0}$ here is reduced to 
$T_{u_b^{-1}}$. Thus the necessary $t$-degree is
$$
\frac{l(w_0)}{2}+(c^+,\rho^\bep)-(l(w_0)-l(u_b))=
l(u_b)-\frac{l(w_0)}{2}+(c^+,\rho^\bep).
$$
We will omit the details (see \cite{CO1}).
We note some connection of the $t$-degrees in 
(\ref{E-ast-rengena}) with the expansion of
the $\mu$-function, the kernel serving the 
nonsymmetric Macdonald polynomials and 
related to the $q$-Kostant function, analyzed 
in certain cases in \cite{Ion2} (see his $D_\lambda$).  

To obtain the second formula in (\ref{ebminusdag}),
we use that 
\begin{align*}
&\ \ \,u_b^{-1}\,\pi_c^{-1}\ =\ w_0\,c^{-1}\ =\ c^\iota\,w_0\ =\ 
\pi_{c^\iota}\,u_{b^\iota}\hbox{\,\, and\, that\, here}\\
&l(u_b^{-1}\,\pi_c^{-1})=
l(u_b^{-1})+l(\pi_c^{-1})=l(\pi_{c^\iota})+l(u_{b^\iota})=
l(\pi_{c^\iota}\,u_{b^\iota}).
\end{align*}
\sq

\subsection{\bf The extremal parts}
Let $\overline{\E}_b$ and $\E^\dag_b$ be the {\em extremal
parts\,} of $\overline{E}_b$ and $E^\dag_b$, i.e. the subsums
of the terms $C_a X_a$, where $C_a\in \Z[q^{\pm 1}]$, for
$a\in W(b)$ such that $a\succeq b$, i.e. those satisfying
(\ref{macdbruhat}).

The polynomial $\overline{\E}_b$ is simply 
$\sum_{a\succeq b} X_a$  
(Proposition 2.4 from \cite{CO1}). It results from 
the formula 
\begin{align*}
&\overline{E}_b=\overline{T}'_{w_0 u_b^{-1}}(E_{b_+}),
\hbox{\, where\, } \overline{T}'_i=\overline{T}_i+1, 
\end{align*}
and $\overline{T}'_{w}$ for $w\in W$ (or any
$\hw\in \hat{W}$) is defined using the 
homogeneous Coxeter relations for $\overline{T}'_i$,
which readily follow from the fact that
$\overline{T}'_i=\lim_{t\to 0} t^{1/2}T_i^{-1}$.
Calculating the extremal parts 
$\overline{\T}_i'$ of $\overline{T}'_i$, one has
\begin{align}
\overline{\E}_b\ =\ \overline{\T}'_{u_b^{-1}w_0}(X_{b_+})\ =\ 
\sum_{a\succeq b} X_a,
\hbox{\, where }&\notag\\
\overline{\T}_i'(X_b)=
\begin{cases}
X_{s_i(b)}+X_b, & \text{if }(b,\al_i)>0,\\
X_b, & \text{if }(b,\al_i)=0,\\
0, & \text{if }(b,\al_i)<0.
\label{intmodsie}
\end{cases}
\end{align}
The calculation of $\overline{\T}'_{u_b^{-1}w_0}(X_{b_+})$
here can be performed by induction or using
the connection between the nil-Hecke algebras
and the Bruhat order. We will need below the following
variant of this calculation.

\begin{lemma}\label{LEMMIX}
For $b\in P_+$ and 
$\,u,v\in W$ such that $\,l(uv)=l(u)+l(v)$,
\begin{align}\label{lemmixed}
\overline{\T}_u\overline{\T}'_v (X_b)&\!=\!
\sum_{w=uv'}X_{w(b)},\hbox{\, where }
l(w)\!=\!l(u)+l(v') \hbox{ and } v\geq v', \\
&\!=\!\!\!\!\!\!\sum_{uv(b)\succeq a\in W(b)}\!X_{a},\hbox{\, where }
(a,\al)<0 \hbox{ for } \{\al\}=\la(u^{-1}),\notag
\end{align}
for the Bruhat order $\geq$ and the ordering $\succeq$ from
(\ref{macdbruhat}).\
\sq 
\end{lemma}
\smallskip

The polynomials $\E^\dag_b$ are significantly more involved
than $\overline{\E}_b$, but they can be linked to the latter 
via the following extension of formula (\ref{rhogb}).

\begin{proposition}\label{EBMEXT}
Let $b\in P$,\ $u=u_{b^\iota}$ and
$v=u_{b^\iota}^{-1}w_0$. 

(i) In the notation of
Proposition \ref{ASTBARDAG} and Lemma \ref{LEMMIX},
\begin{align}
\E_b^\dag\ &=\ q^{-(c_+,b_+)}\,
\Bigl(\,\pi_{r^\iota}\ \overline{\T}_{j_l^\iota}^{\,\ep_l}\cdots
\,\overline{\T}_{j_2^\iota}^{\,\ep_2}
\overline{\T}_{j_1^\iota}^{\,\ep_1}\,
\bigl(\overline{\T}_u
\overline{\T}_v^\prime(X_{b^\iota_+})\bigr)\Bigr)^\ast,
\label{ebminusext}
\end{align}
where $\overline{\T}_u\overline{\T}_v^\prime(X_{b^\iota_+})=
\sum_w\,X_{w^\iota(b^\iota_+)}\,$ for the summation over $\,w\in W\,$ 
such that $\,l(w)=l(u_b)+l(u_b^{-1}\,\,w)$; equivalently,
$\overline{\T}_u\overline{\T}_v^\prime(X_{b^\iota_+})=
\sum_{a\in W(b^\iota)}X_{a}$ provided $(a,\al)<0$
for $\al\in \la(u_{b^\iota}^{-1})$. 

(ii) Taking here $\,c=u_{b^\iota}^{-1}(\check{\rho}^\bep)\,$,
one has 
\begin{align}\label{picrho}
&\pi_c=u_{b^\iota}^{-1}\pi_{\check{\rho}^\bep},\hbox{ where }  
l(\pi_c)=l(u_{b^\iota})+l(\pi_{\check{\rho}^\bep}),\notag\\
&\E_b^\dag= q^{-(\check{\rho}^\bep,b_+)}\,
\Bigl(\,\overline{\T}_{u_{b^\iota}^{-1}}\ 
\overline{\T}^\prime_{\pi_{\check{\rho}^\bep}}\,
\bigl(\,\overline{\T}_u
\overline{\T}_v^\prime(X_{b^\iota_+})\,\bigr)\Bigr)^\ast.
\end{align} 
\sq
\end{proposition}

If $\ep_p=\varnothing\,$ for some $\,p\,$ in
(\ref{ebminusext}), then there can 
be generally negative terms upon applying the corresponding 
$\overline{\T}_{j_p^\iota}$;
subtract $X_b$ from (\ref{intmodsie}) to obtain the formula
for $\overline{T_i}$ (without $\,\prime$).
Therefore the argument that proved Theorem \ref{EXTDAG} can not 
be immediately used to establish that all nonzero coefficients 
of $\E_b^\dag$ are in the form $q^{-m}$ for $m\ge 0$. 

This actually holds true for any $b\in P$,
but the justification is more involved; it is based 
on the method used in the proof of Theorem \ref{EXTDAG} 
and formula (2.52) from \cite{CO1} and will be not  
discussed here. However the following application of 
(\ref{ebminusext}) is straightforward.

\begin{proposition}\label{PUREQPAR}
Let $u_b\,$ be $\,w^I_0\,$ for
the root subsystem corresponding to a subset $I\subset \Ga$,
including the case $I=\emptyset$. Then $u_c=w_0w^I_0$\,
for $c=\sum_{i\not\in I}\om_i^\vee$ for $\,\bep\,=\varnothing$
and  $c=\sum_{i\not\in I}\om_i$ in the twisted case;
cf. (\ref{xwo})). Formula (\ref{ebminusext}) holds
for $b=w^I_0(b_-)$ and such $\,c\,$; moreover, 
$\ep_p=\prime\,$ for all indices $p$.  In this case,
$u^\iota=w_0^I,$ \,$v^\iota=w_0^I w_0\,$ and
\begin{align*}
\overline{\T}_u\overline{\T}_v^\prime(X_{b^\iota_+})= 
\sum_w\,X_{w^\iota(b^\iota_+)}\,,\hbox{\, where the summation
is over\, } w\in W
\end{align*}
such that $\,l(w)\!=\!l(w^I_0)\!+\!l(w^I_0\,w).$

Thus the nonzero coefficients of $X_{c}$ in
$\E^\dag_b$ are all in the form $q^{-e(-b,w)}$, where $c=w(b)$ 
and $e(-b,w)\in \Z_+$.
Moreover, 
$e(-a - b,w)\!=\!e(-a,w)+e(-b,w)$ provided  
$\,u_{a}\!=\! w^I_0=u_b$.
\sq
\end{proposition}

The simplest case of the proposition is when  $I=\Ga$
and $u_b=w_0$.
Then $c=0\,$ and $u_c=$id$=\pi_c$; we obtain that $\E^\dag_b=X_b$
for regular {\em dominant} $b$, which of course holds for 
any $b\in P_+$ (not only regular) due to (\ref{macdbruhat}).
 
If $I=\emptyset$ then $u_b=$id
and $b\in P_-$. This case is covered by 
Theorems \ref{EXTDAG},\ref{EXTADDIT}; here $c=\check{\rho}^\bep$.
We note that the monomials $X_a$ for
$a\succeq b$ can be missing (with zero coefficients) in $\E^\dag_b$
if $b\not\in P_-$.
\medskip

{\em On embeddings of dag-polynomials.}
Without going into detail, let us mention another
special case of Proposition \ref{PUREQPAR} 
directly related to Theorem \ref{ADDWITHIN}. 
The decomposition of $\E_{b_-}^\dag$ in terms of 
the blocks from $\mathfrak{T}^{\bep}(X_{a_\circ})$ 
is directly related with $\E_{b}^\dag$ for $b=w_0^I(b_-)$ 
when the sets $I\subset \Ga$ are 
unions of disconnected points (i.e. contain no nontrivial
segments).  

Such polynomials $\E_{b}^\dag$ are naturally embedded 
into  $\E_{b_-}^\dag$ upon the multiplication 
by $\,q^{m}$ for $m=\sum_{i\in I}\,(b_-\,,\check{\al}^\bep_i)\,$,
where $\check{\al}^\bep_i$ is $\al_i^\vee$ in the untwisted case
and $\al_i$ otherwise.
Then they become sums of the blocks described in
Theorem \ref{ADDWITHIN} for $u=\prod_{i\in I'}\,s_i\,$
over all the subsets $I'\supseteq I$ (without nontrivial segments).
The images of $\E_{a}^\dag$  and $\E_{b}^\dag$  in
$\E_{b_-}^\dag$ associated 
with $I_a$ and $I_b$ ``intersect" exactly by the image of 
the $E$-dag polynomial corresponding to $I_a\cup I_b$ (if the latter 
set contains no segments). 
\smallskip

The best way to justify this is by using (\ref{picrho}) for
$u_{b^\iota}=\prod_{i\in I}\,s_i\,$. Let us take one
$s_{i_\circ}$ here for the sake of definiteness; then
$b=s_i(b_-)$, $\,pi_c=s^\iota_{i}\pi_{\check{\rho}^\bep}$.
Setting
$$\Si=X_{b_-}\sum_{I\subseteq \Ga} \prod_{i\in I} X_{\al_i}
\hbox{\  \,and\, \ }
\Si^\circ=X_{b_-}\!\!\sum_{i_\circ\in
I\subseteq \Ga} \prod_{i\in I} X_{\al_i},
$$
formula (\ref{picrho}) and  Theorem \ref{ADDWITHIN}
result in
\begin{align}\label{picrhone}
&\E_b^\dag= q^{-(\check{\rho}^\bep,b_+)}\,
\Bigl(\,\overline{\T}_{i_\circ^\iota}\ 
\overline{\T}^\prime_{\pi_{\check{\rho}^\bep}}\,
\bigl(\,\Si-\Si^{\circ}\,\bigr)\Bigr)^\ast,
\end{align} 
Using formula (2.52) from \cite{CO1} (and notation
$(\overline{T}_{i}^\dag)'$ there),
\begin{align*}
\Bigl(\,\overline{\T}_{i_\circ^\iota}\ 
\overline{\T}^\prime_{\pi_{\check{\rho}^\bep}}\,
\bigl(\,\Si\,\bigr)\Bigr)^\ast=
(\overline{\T}_{i_\circ^\iota}^\dag)'\,(\E_{b_-}^\dag)=
(1-q^{(b,\check{\al}_{i_\circ}^\bep)})\,\E_{b}^\dag.
\end{align*} 
Then $\overline{\T}_{i_\circ^\iota}\ 
\overline{\T}^\prime_{\pi_{\check{\rho}^\bep}}\,
\bigl(\,\Si^{\circ}\,\bigr)=
 -\overline{\T}^\prime_{\pi_{\check{\rho}^\bep}}\,
\bigl(\,\Si^{\circ}\,\bigr),$\ which gives that
$$
q^{-(\check{\rho}^\bep,b_+)}\,
\Bigl(\,(\overline{\T}^\prime_{\pi_{\check{\rho}^\bep}}\,
\bigl(\,\Si^{\circ}\,\bigr)\Bigr)^\ast=
q^{(b,\check{\al}_{i_\circ}^\bep)}\,\E_b^\dag,
$$
where the left-hand side is \,
$q^{(\check{\rho}^\bep,b_-)}\,\mathfrak{T}^\bep\,(X_{s_i(b)})$\,
from 
Theorem \ref{ADDWITHIN} in the case of $u=s_{i_\circ}$.
\smallskip

We do not have any conjectures so far concerning the
representation-theoretic meaning of $E_b^\dag$ or 
$\E_b^\dag$ apart from antidominant $b$, thought 
all coefficients of full $E_b^\dag$  are expected to be positive 
for any $b\in P$, which was conjectured in \cite{CO1}.
\medskip

\subsection{\bf The affine conjecture} Let $\hat{\mathfrak{g}}$
be the Kac-Moody algebra associated with $\tR^\bep$,
$\hat{\mathfrak{b}}_{+}$ and $\hat{\mathfrak{n}}_{+}$ its
Borel subalgebra and its nilpotent subalgebra.
Following Conjecture 2.7 from \cite{CO1}, we
will consider only the simply-laced cased. Thus 
$\bep\,=\varnothing\,$
and $\hat{\mathfrak{g}}\!=\!\mathfrak{g}[z,z^{-1}]\oplus \C c$\,
with the standard central element $c$ and the commutator;\,
$\mathfrak{g}$ is the simple Lie algebra  associated with $R$.

Given $b\in P_-$, let $\la_\circ$ be a minuscule weight
$\om_r$ (for $r\in O'$) or zero such that $a=\la_\circ-b\in Q_+$. 
Consider the level-one irreducible integrable  
$\hat{\mathfrak{g}}$-module $L_{\tilde{\la}_\circ}$ 
for $\tilde{\la}_\circ=[\la_\circ,1]$. It is generated by
the vacuum vector $v_\circ$ defined with respect to 
$\hat{\mathfrak{b}}_{+}$, i.e. satisfying 
$\hat{\mathfrak{n}}_+(v_\circ)=0$.
\smallskip

{\em The Demazure module} $\d_b$ is the
following $\,\hat{\mathfrak{b}}_{+}$-module:
$$
U(\hat{\mathfrak{b}}_{+})(v_{b})\subset L_{\tilde{\la}_\circ},
\hbox{\, where \,} b=\la_\circ-a=(-a)(\la_\circ);
$$ 
we set $v_b=(-a)(v_{\circ})$. More generally, 
$\la=\tilde{w}(\la_\circ)$ and $ v_\la=\tilde{w}(v_{\la_\circ})$
for the standard action of $\tilde{w}\in \tilde{W}$ in $P$
(with $Q\subset \tilde{W}$ acting
via translations) and its lift to $L_{\tilde{\la}_\circ}.$  
Due to considering only antidominant $b$, the Demazure module
$\d_b$ is $\mathfrak{g}$-invariant.

This module has the {\em standard Kac-Moody 
grading $d(v)$}. This degree will be counted from $v_b$;
for instance, $\ d(e_{\al_0}^m v_b)=m\,$ for $m\ge 0$.
The {\em PBW-degree\,} of $\,v\in \d_b\,$  is defined 
as the minimal number $\de=\de(v)$ such that 
$v\in U(\hat{\mathfrak{b}}_{+})_\de\,(v_{b})$, where 
$U(\hat{\mathfrak{b}}_{+})_\de$ is the $\de$-th piece of the
PBW filtration on $U(\hat{\mathfrak{b}}_{+})$.
\smallskip

For any $c\in Q_+$, let $C^\dag_c(b)$ be
the $q$-character of the weight space
$\,\d_b(b+c)=
L_{\tilde{\la}_\circ}(b+c)\cap \d_b\,$ with respect to
the composite degree $\,-d(v)-\de(v)\,$ (note the minuses).
To be more exact, we consider here the graded module
of $\d_b(b+c)$ with respect to the Kac-Moody- and 
PBW-filtrations and its character.

\begin{conjecture}\label{CONJDEM}
For $\,b\in P_-\,, $ 
\begin{align}\label{conjdema}
\,\ \ \ \ \ \ \ E^\dag_{b}\ =\ 
\sum_{c\in Q_+} C^\dag_c(b)X_{b+c}\,,
\hbox{\, \,where\, \,} b+c\,\succeq\, b.\ \ \ \ \ \ \ 
\hbox{\sq}
\end{align}
\end{conjecture}

Defining the characters $\overline{C}_c(b)$ only for $-d(v)$
(without $\de(v)$), 
the following formula is due to \cite{San,Ion1}:
\begin{align}\label{conjdemba}
&\overline{E}_b\,=\,\sum_{c\in Q_+} \overline{C}_c(b)X_{b+c} 
\hbox{\, \ for\, \ } b\in P_-,\ \,b+c\,\succeq\, b.
\end{align}
As a matter of fact, $b$ can be taken in (\ref{conjdemba}) 
from $P$ and any twisted root systems is allowed here.
\smallskip

Restricting ourselves to the $W$-extremal vectors (their  
$d$-grading is zero), the conjecture states that
$q^{-d(-b,w)}=C^\dag_{w(b)-b}=q^{-e(-b,w)}$, which was proven
above for all classical root systems in the untwisted setting
(including $B,C$) and $G_2$.
We do not suggest any affine conjecture(s) apart from the
$A\!D\!E$-systems in this paper,
though the theory of PBW-filtration seems quite doable in 
the twisted setting.

\medskip

\subsection{\bf The systems \texorpdfstring{$A_1-A_3$}{A1-A3}} 
We generally do not have 
systematic tools for calculating the PBW-filtration in
the Kac-Moody case. However this can be done
in sufficiently small examples. We note that
Conjecture \ref{CONJDEM} can be checked partially,
since it implies certain inequalities (sometimes
strong enough) for the $q$-degrees of
$C^\dag_c(b)$ vs. $\overline{C}_c(b)$. These inequalities were
checked numerically in quite a few cases.

The examples provided below are of theoretical nature;
the PBW-calculations were performed using the explicit 
realization of the level one Demazure modules 
(see \cite{CL}, \cite{FoL}).
\smallskip

{\em The case of $A_1$.} Then the $E$-polynomials
$E^\dag$-polynomial can be computed explicitly (see e.g. 
formula (1.30) from \cite{CO2}):
\[
E^\dag_{-n}(q,X)=\sum_{j=0}^n X^{2j-n}q^{-j}\binom{n}{j}_{q^{-1}},\ 
\binom{n}{j}_q\equal\frac{(1-q^{n-j+1})
\dots(1-q^n)}{(1-q)
\dots(1-q^j)}.
\]
This proves Conjecture \ref{CONJDEM} in type $A_1$ (see e.g. \cite{CL} 
for the character of the corresponding Demazure module).
Indeed, the coefficient of $\,X^k\,$ corresponds to the character of the
$\msl_2$-weight subspace of weight $\,k$ of the Demazure 
module. The weight of the cyclic vector is exactly $\,-n\,$,
which corresponds to the term $X^{-n}$ at  $j=0$. It is easy to see 
that the PBW-degrees of all vectors of the $\msl_2$-weight 
equal to $(-n+2j)$ are $\,j\,$. Therefore it suffices to note that 
the $q$-binomial coefficient $\binom{n}{j}_{q^{-1}}$ is exactly 
the Kac-Moody character of the subspace of weight $\,j$.
\smallskip

{\em An example for $A_2$}. 
Let $\la_\circ=0$ and $b=-2\th=-2\om_1-2\om_2.$
I.e. $\d_{-2\th}$ is inside
$L_{\tilde{\la}_\circ= [0,1]}$\,.
We will discuss in this example only the zero-level
subspace of $L_{[0,1]}$, the eigenspace through $v_{\la_\circ}$; 
thus $c=2\th$. 

The zero-level subspace of the Demazure module
$\d_b$ has the following $d$-character (that for the grading $d(v)$;
note the sign): 
\begin{align}\label{barA2}
\,\,3\, +\, 4q\, +\, 5q^2\, +\, 2q^3\, +\, q^4\ \, =\, 
\hbox{CT}(\overline{E}_{-2\th})\mid_{q\mapsto q^{-1}}\,,
\end{align}
where CT$=$Constant Term.
The character for $d(v)+\de(v)$ equals
\begin{align}\label{dagA2}
q^2 + 2q^3 + 6q^4 + 4q^5 + 2q^6\ =\ 
\hbox{CT}(E^\dag_{-2\th})\,\mid_{q\mapsto q^{-1}}\,.
\end{align}
Let us provide the {\em full bi-character} at zero-level
of $\d_{-2\th}$ corresponding to the weight function
$\,q^{d(v)}\,\tau^{\de(v)}\,$:
$$
\tau^2\bigl((1+\tau+\tau^2)  + (1+3\tau)q   +  (2+3\tau)q^2 +   
(1+\tau)q^3 +q^4\bigr).
$$
It does coincide with (\ref{dagA2}) when $\tau=q$. We do not have
any general conjectures concerning the full bi-characters of
level-one Demazure modules.
\medskip

{\em Examples for $A_3$}.
Let us list the full bi-characters (for $\,q^{d(v)}\,
\tau^{\de(v)}\,$) and the corresponding $E$-polynomials
for small $b$ in the case of $A_3$, namely, for $b=-2\om_1$,
$b=-\om_1-\om_3$ and $b=-\om_1-\om_2$.
We set $A[l,m,n]=X_{\al_1}^l X_{\al_2}^m X_{\al_3}^n$.
\smallskip

{\em The case of $b=-2\om_1$.}
The full bi-character equals

\renewcommand{\baselinestretch}{0.5} 
{\small
\noindent
\(
X_{b}\Bigl(A[0,0,0]+(\tau+q \tau) A[1,0,0]+(\tau+q \tau) A[1,1,0]
+(\tau+q \tau) A[1,1,1]\\+\tau^2 A[2,0,0]+\left(\tau^2+
q \tau^2\right) A[2,1,0]
+\left(\tau^2+q \tau^2\right) A[2,1,1]+\tau^2 A[2,2,0]\\+
\left(\tau^2+q \tau^2\right) A[2,2,1]+\tau^2 A[2,2,2]\Bigl).
\)
}
\renewcommand{\baselinestretch}{1.2} 

Upon the substitution $\tau\mapsto q\,$, the bi-character
becomes $\E_b^\dag\mid_{q\mapsto q^{-1}}\,$ here
and below.
Let us provide the corresponding $E$-polynomial (we use SAGE
software for this and the next two $E$-polynomials):

$E_{-2\om_1}(X;q,t)\!=\!$
\renewcommand{\baselinestretch}{0.5} 
{\small
\noindent
\(
X_{-2\om_1}\Bigl(
A[0,0,0]+\frac{(1+q) (1-t) A[1,0,0]}{1-q^2 t}+\frac{(1+q) (1-t) 
A[1,1,0]}{1-q^2 t}\\
+\frac{(1+q) (1-t) A[1,1,1]}{1-q^2 t}
+\frac{(1-t) A[2,0,0]}{1-q^2 t}+\frac{(1+q) (1-t)^2 
A[2,1,0]}{(1-q t) \left(1-q^2 t\right)}+\frac{(1+q) (1-t)^2 
A[2,1,1]}{(1-q t) \left(1-q^2 t\right)}\\+\frac{(1-t) 
A[2,2,0]}{1-q^2 t}
+\frac{(1+q) (1-t)^2 A[2,2,1]}{(1-q t) 
\left(1-q^2 t\right)}+\frac{(1-t) A[2,2,2]}{1-q^2 t}\Bigr).
\)
}
\renewcommand{\baselinestretch}{1.2} 

Setting here and below 
$t\to\infty$ and $q\mapsto q^{-1}$, one obtains
the bi-character where $\tau=q$.

\smallskip

{\em The case of $b=-\om_1-\om_3$.}
The full bi-character equals

\renewcommand{\baselinestretch}{0.5} 
{\small
\noindent
\(
X_b\Bigl(
A[0,0,0]+\tau A[0,0,1]+\tau A[0,1,1]+\tau A[1,0,0]
+\tau^2 A[1,0,1]+\tau A[1,1,0]\\+\tau (1+q+2 \tau) A[1,1,1]
+\tau^2 A[1,1,2]+\tau^2 A[1,2,1]+\tau^2 A[1,2,2]
\\+\tau^2 A[2,1,1]+\tau^2 A[2,2,1]+\tau^2 A[2,2,2]
\Bigl);
\)
}
\renewcommand{\baselinestretch}{1.2}

$X_{\om_1+\om_3}E_{-\om_1-\om_3}(X;q,t)\!=\!$
\renewcommand{\baselinestretch}{0.5} 
{\small
\noindent
\(
A[0,0,0]+\frac{(1-t) \bigl(A[0,0,1]+A[0,1,1]+A[1,0,0]\bigr)
}{1-q t}
\\+\frac{(1-t)^2 A[1,0,1]}{(1-q t)^2}+\frac{(1-t) 
A[1,1,0]}{1-q t}+\frac{(1-t) \left(3+q-2 t-2 q t-4 q^2 t^3
+3 q^2 t^4+q^3 t^4\right) A[1,1,1]}{(1-q t)^2 
\left(1-q^2 t^3\right)}\\
+\frac{(1-t)^2 A[1,1,2]}{(1-q t)^2}
+\frac{(1-t)^2 A[1,2,1]}{(1-q t)^2}
+\frac{(1-t)^2 A[1,2,2]}{(1-q t)^2}+\frac{(1-t)^2 
A[2,1,1]}{(1-q t)^2}+\frac{(1-t)^2 A[2,2,1]}{(1-q t)^2}
\\+\frac{(1-t) \left(1-t+t^3-2 q t^3+q^2 t^4\right) 
A[2,2,2]}{(1-q t)^2 \left(1-q^2 t^3\right)}.
\)
}
\renewcommand{\baselinestretch}{1.2} 
\smallskip

{\em The case of $b=-\om_1-\om_2$.}
The full bi-character equals

\renewcommand{\baselinestretch}{0.5} 
{\small
\noindent
\(
X_b\Bigl(
A[0,0,0]+\tau A[0,1,0]+\tau A[0,1,1]+\tau A[1,0,0]+
\tau (1+q+\tau) A[1,1,0]
\\+\tau (1+q+\tau) A[1,1,1]+\tau^2 A[1,2,0]+(2+q) \tau^2 A[1,2,1]
+\tau^2 A[1,2,2]
\\+\tau^2 A[2,1,0]
+\tau^2 A[2,1,1]+\tau^2 A[2,2,0]+\tau^2 (1+q+\tau) A[2,2,1]
+\tau^2 A[2,2,2]\\+\tau^3 A[2,3,1]+\tau^3 A[2,3,2]
\Bigl);
\)
}
\renewcommand{\baselinestretch}{1.2}

$X_{\om_1+\om_2}E_{-\om_1-\om_2}(X;q,t)\!=\!$
\renewcommand{\baselinestretch}{0.5} 
{\small
\noindent
\(
A[0,0,0]+\frac{(1-t) \bigl(A[0,1,0]+A[0,1,1]+A[1,0,0]\bigr)
}{1-q t}
\\+\frac{(1-t) \left(2+q-t-2 q t-3 q^2 t^2+2 q^2 t^3+q^3 t^3\right) 
\bigl(A[1,1,0]+A[1,1,1]\bigr)}{(1-q t)^3 \left(1+q t\right)}
+\frac{(1-t)^2 
A[1,2,0]}{(1-q t)^2}\\+\frac{(1-t)^2 \left(2+q-q t^2
-2 q^2 t^2\right) A[1,2,1]}{(1-q t)^3 (1+q t)}
+\frac{(1-t)^2 A[1,2,2]}{(1-q t)^2}+\frac{(1-t)^2 
A[2,1,0]}{(1-q t)^2}+\frac{(1-t)^2 A[2,1,1]}{(1-q t)^2}
\\+\frac{(1-t) \left(1-t+t^2-2 q t^2+q^2 t^3\right) 
A[2,2,0]}{(1-q t)^3 (1+q t)}+\frac{(1-t)^2 \left(2+q-t+q t
-2 q t^2-q^2 t^2\right) A[2,2,1]}{(1-q t)^3 (1+q t)}
\\+\frac{(1-t) \left(1-t+t^2-2 q t^2+q^2 t^3\right) 
A[2,2,2]}{(1-q t)^3 (1+q t)}+\frac{(1-t)^2 \left(1-q t^2\right) 
A[2,3,1]}{(1-q t)^3 (1+q t)}+\frac{(1-t)^2 
\left(1-q t^2\right) A[2,3,2]}{(1-q t)^3 (1+q t)}.
\)
}
\renewcommand{\baselinestretch}{1.2} 
\medskip

\section{\sc Appendix: the system \texorpdfstring{$F_4$}{F4}}

Let us provide formulas for the dag-polynomials in the untwisted
and twisted cases of $F_4$. Recall (\ref{Pomwom}) and that 
$$
A[c]=\prod_{i=1}^n \,X_{\al_i}^{\,c_i} \hbox{\, for\, }
Q\ni c=\sum_{i=1}\, c_i \al_i.
$$
\smallskip

{\em Twisted $\tilde{\E}^\dag$-polynomials for $F_4$.}
\medskip

$(\tilde{\E}^\dag_1)^\nu=$
\renewcommand{\baselinestretch}{0.5} 
{\small
\noindent
\(
1
+\frac{A[1,0,0,0]}{q^2}+\frac{A[1,1,0,0]}{q^2}+\frac{A[1,1,2,0]}{q^2}
+\frac{A[1,1,2,2]}{q^2}+\frac{A[1,2,2,0]}{q^2}+\frac{A[1,2,2,2]}{q^2}
+\frac{A[1,2,4,2]}{q^2}+\frac{A[1,3,4,2]}{q^2}+\frac{A[2,2,2,0]}{q^2}
+\frac{A[2,2,2,2]}{q^2}+\frac{A[2,2,4,2]}{q^2}+\frac{A[2,4,4,2]}{q^2}
+\frac{A[2,4,6,2]}{q^2}+\frac{A[2,4,6,4]}{q^2}+\frac{A[3,3,4,2]}{q^4}
+\frac{A[3,4,4,2]}{q^4}+\frac{A[3,4,6,2]}{q^4}+\frac{A[3,4,6,4]}{q^4}
+\frac{A[3,5,6,2]}{q^4}+\frac{A[3,5,6,4]}{q^4}
+\frac{A[3,5,8,4]}{q^4}+\frac{A[3,6,8,4]}{q^4}+\frac{A[4,6,8,4]}{q^4}.
\)
}
\renewcommand{\baselinestretch}{1.2} 
\smallskip

$(\tilde{\E}^\dag_2)^\nu=$
\renewcommand{\baselinestretch}{0.5} 
{\small
\noindent
\(
1+\frac{A[0,1,0,0]}{q^2}+\frac{A[0,1,2,0]}{q^2}+\frac{A[0,1,2,2]}{q^2}
+\frac{A[0,2,2,0]}{q^2}+\frac{A[0,2,2,2]}{q^2}+\frac{A[0,2,4,2]}{q^2}
+\frac{A[0,3,4,2]}{q^4}+\frac{A[1,1,0,0]}{q^2}+\frac{A[1,1,2,0]}{q^2}
+\frac{A[1,1,2,2]}{q^2}+\frac{A[1,3,2,0]}{q^4}+\frac{A[1,3,2,2]}{q^4}
+\frac{A[1,3,4,0]}{q^4}+\frac{A[1,3,4,4]}{q^4}+\frac{A[1,3,6,2]}{q^4}
+\frac{A[1,3,6,4]}{q^4}+\frac{A[1,5,6,2]}{q^4}+\frac{A[1,5,6,4]}{q^4}
+\frac{A[1,5,8,4]}{q^4}+\frac{A[2,2,2,0]}{q^2}+\frac{A[2,2,2,2]}{q^2}
+\frac{A[2,2,4,2]}{q^2}+\frac{A[2,3,2,0]}{q^4}+\frac{A[2,3,2,2]}{q^4}
+\frac{A[2,3,4,0]}{q^4}+\frac{A[2,3,4,4]}{q^4}+\frac{A[2,3,6,2]}{q^4}
+\frac{A[2,3,6,4]}{q^4}+\frac{A[2,4,4,0]}{q^4}+\frac{A[2,4,4,4]}{q^4}
+\frac{A[2,4,8,4]}{q^4}+\frac{A[2,5,4,2]}{q^4}+\frac{A[2,5,8,2]}{q^4}
+\frac{A[2,5,8,6]}{q^4}+\frac{A[2,6,6,2]}{q^4}+\frac{A[2,6,6,4]}{q^4}
+\frac{A[2,6,8,2]}{q^4}+\frac{A[2,6,8,6]}{q^4}+\frac{A[2,6,10,4]}{q^4}
+\frac{A[2,6,10,6]}{q^4}+\frac{A[2,7,8,4]}{q^6}
+\frac{A[2,7,10,4]}{q^6}
+\frac{A[2,7,10,6]}{q^6}+\frac{A[3,3,4,2]}{q^4}+\frac{A[3,5,4,2]}{q^4}
+\frac{A[3,5,8,2]}{q^4}+\frac{A[3,5,8,6]}{q^4}+\frac{A[3,7,8,2]}{q^6}
+\frac{A[3,7,8,6]}{q^6}+\frac{A[3,7,12,6]}{q^6}
+\frac{A[3,9,12,6]}{q^6}
+\frac{A[4,5,6,2]}{q^4}+\frac{A[4,5,6,4]}{q^4}+\frac{A[4,5,8,4]}{q^4}
+\frac{A[4,6,6,2]}{q^4}+\frac{A[4,6,6,4]}{q^4}+\frac{A[4,6,8,2]}{q^4}
+\frac{A[4,6,8,6]}{q^4}+\frac{A[4,6,10,4]}{q^4}
+\frac{A[4,6,10,6]}{q^4}
+\frac{A[4,7,8,2]}{q^6}+\frac{A[4,7,8,6]}{q^6}+\frac{A[4,7,12,6]}{q^6}
+\frac{A[4,8,8,4]}{q^4}+\frac{A[4,8,12,4]}{q^4}
+\frac{A[4,8,12,8]}{q^4}
+\frac{A[4,9,10,4]}{q^6}+\frac{A[4,9,10,6]}{q^6}
+\frac{A[4,9,12,4]}{q^6}
+\frac{A[4,9,12,8]}{q^6}+\frac{A[4,9,14,6]}{q^6}
+\frac{A[4,9,14,8]}{q^6}
+\frac{A[4,10,12,6]}{q^6}+\frac{A[4,10,14,6]}{q^6}
+\frac{A[4,10,14,8]}{q^6}+\frac{A[5,7,8,4]}{q^6}
+\frac{A[5,7,10,4]}{q^6}
+\frac{A[5,7,10,6]}{q^6}+\frac{A[5,9,10,4]}{q^6}
+\frac{A[5,9,10,6]}{q^6}
+\frac{A[5,9,12,4]}{q^6}+\frac{A[5,9,12,8]}{q^6}
+\frac{A[5,9,14,6]}{q^6}
+\frac{A[5,9,14,8]}{q^6}+\frac{A[5,11,14,6]}{q^8}
+\frac{A[5,11,14,8]}{q^8}+\frac{A[5,11,16,8]}{q^8}
+\frac{A[6,9,12,6]}{q^6}+\frac{A[6,10,12,6]}{q^6}
+\frac{A[6,10,14,6]}{q^6}+\frac{A[6,10,14,8]}{q^6}
+\frac{A[6,11,14,6]}{q^8}+\frac{A[6,11,14,8]}{q^8}
+\frac{A[6,11,16,8]}{q^8}+\frac{A[6,12,16,8]}{q^8}.
\)
}
\renewcommand{\baselinestretch}{1.2} 
\smallskip

$(\tilde{\E}^\dag_3)^\nu=$
\renewcommand{\baselinestretch}{0.5} 
{\small
\noindent
\(
1+\frac{A[0,0,1,0]}{q}+\frac{A[0,0,1,1]}{q}+\frac{A[0,1,1,0]}{q}
+\frac{A[0,1,1,1]}{q}+\frac{A[0,1,2,0]}{q^2}+\frac{A[0,1,2,2]}{q^2}
+\frac{A[0,1,3,1]}{q^2}+\frac{A[0,1,3,2]}{q^2}+\frac{A[0,2,3,1]}{q^2}
+\frac{A[0,2,3,2]}{q^2}+\frac{A[0,2,4,2]}{q^2}+\frac{A[1,1,1,0]}{q}
+\frac{A[1,1,1,1]}{q}+\frac{A[1,1,2,0]}{q^2}+\frac{A[1,1,2,2]}{q^2}
+\frac{A[1,1,3,1]}{q^2}+\frac{A[1,1,3,2]}{q^2}+\frac{A[1,2,2,0]}{q^2}
+\frac{A[1,2,2,2]}{q^2}+\frac{A[1,2,3,0]}{q^3}+\frac{A[1,2,3,3]}{q^3}
+\frac{A[1,2,5,2]}{q^3}+\frac{A[1,2,5,3]}{q^3}+\frac{A[1,3,3,1]}{q^2}
+\frac{A[1,3,3,2]}{q^2}+\frac{A[1,3,5,1]}{q^3}+\frac{A[1,3,5,4]}{q^3}
+\frac{A[1,3,6,2]}{q^3}+\frac{A[1,3,6,4]}{q^3}+\frac{A[1,4,5,2]}{q^3}
+\frac{A[1,4,5,3]}{q^3}+\frac{A[1,4,6,2]}{q^3}+\frac{A[1,4,6,4]}{q^3}
+\frac{A[1,4,7,3]}{q^3}+\frac{A[1,4,7,4]}{q^3}+\frac{A[2,2,3,1]}{q^2}
+\frac{A[2,2,3,2]}{q^2}+\frac{A[2,2,4,2]}{q^2}+\frac{A[2,3,3,1]}{q^2}
+\frac{A[2,3,3,2]}{q^2}+\frac{A[2,3,5,1]}{q^3}+\frac{A[2,3,5,4]}{q^3}
+\frac{A[2,3,6,2]}{q^3}+\frac{A[2,3,6,4]}{q^3}+\frac{A[2,4,4,2]}{q^2}
+\frac{A[2,4,5,1]}{q^3}+\frac{A[2,4,5,4]}{q^3}+\frac{A[2,4,7,2]}{q^3}
+\frac{A[2,4,7,5]}{q^3}+\frac{A[2,4,8,4]}{q^4}+\frac{A[2,5,6,2]}{q^3}
+\frac{A[2,5,6,4]}{q^3}+\frac{A[2,5,7,2]}{q^3}+\frac{A[2,5,7,5]}{q^3}
+\frac{A[2,5,9,4]}{q^4}+\frac{A[2,5,9,5]}{q^4}+\frac{A[2,6,8,4]}{q^4}
+\frac{A[2,6,9,4]}{q^4}+\frac{A[2,6,9,5]}{q^4}+\frac{A[3,4,5,2]}{q^3}
+\frac{A[3,4,5,3]}{q^3}+\frac{A[3,4,6,2]}{q^3}+\frac{A[3,4,6,4]}{q^3}
+\frac{A[3,4,7,3]}{q^3}+\frac{A[3,4,7,4]}{q^3}+\frac{A[3,5,6,2]}{q^3}
+\frac{A[3,5,6,4]}{q^3}+\frac{A[3,5,7,2]}{q^3}+\frac{A[3,5,7,5]}{q^3}
+\frac{A[3,5,9,4]}{q^4}+\frac{A[3,5,9,5]}{q^4}+\frac{A[3,6,7,3]}{q^3}
+\frac{A[3,6,7,4]}{q^3}+\frac{A[3,6,9,3]}{q^3}+\frac{A[3,6,9,6]}{q^3}
+\frac{A[3,6,10,4]}{q^4}+\frac{A[3,6,10,6]}{q^4}
+\frac{A[3,7,9,4]}{q^4}
+\frac{A[3,7,9,5]}{q^4}+\frac{A[3,7,10,4]}{q^4}
+\frac{A[3,7,10,6]}{q^4}
+\frac{A[3,7,11,5]}{q^5}+\frac{A[3,7,11,6]}{q^5}
+\frac{A[4,6,8,4]}{q^4}
+\frac{A[4,6,9,4]}{q^4}+\frac{A[4,6,9,5]}{q^4}+\frac{A[4,7,9,4]}{q^4}
+\frac{A[4,7,9,5]}{q^4}+\frac{A[4,7,10,4]}{q^4}
+\frac{A[4,7,10,6]}{q^4}
+\frac{A[4,7,11,5]}{q^5}+\frac{A[4,7,11,6]}{q^5}
+\frac{A[4,8,11,5]}{q^5}
+\frac{A[4,8,11,6]}{q^5}+\frac{A[4,8,12,6]}{q^6}.
\)
}
\renewcommand{\baselinestretch}{1.2} 
\smallskip

$(\tilde{\E}^\dag_4)^\nu=$
\renewcommand{\baselinestretch}{0.5} 
{\small
\noindent
\(
1+\frac{A[0,0,0,1]}{q}+\frac{A[0,0,1,1]}{q}+\frac{A[0,1,1,1]}{q}
+\frac{A[0,1,2,1]}{q}+\frac{A[0,1,2,2]}{q^2}+\frac{A[1,1,1,1]}{q}
+\frac{A[1,1,2,1]}{q}+\frac{A[1,1,2,2]}{q^2}+\frac{A[1,2,2,1]}{q}
+\frac{A[1,2,2,2]}{q^2}+\frac{A[1,2,3,1]}{q}+\frac{A[1,2,3,3]}{q^2}
+\frac{A[1,2,4,2]}{q^2}+\frac{A[1,2,4,3]}{q^2}+\frac{A[1,3,4,2]}{q^2}
+\frac{A[1,3,4,3]}{q^2}+\frac{A[1,3,5,3]}{q^2}+\frac{A[2,3,4,2]}{q^2}
+\frac{A[2,3,4,3]}{q^2}+\frac{A[2,3,5,3]}{q^2}+\frac{A[2,4,5,3]}{q^2}
+\frac{A[2,4,6,3]}{q^2}+\frac{A[2,4,6,4]}{q^2}.
\)
}
\renewcommand{\baselinestretch}{1.2} 
\medskip

{\em Untwisted $F_4$.}
\smallskip

$\tilde{\E}^\dag_1=$
\renewcommand{\baselinestretch}{0.5} 
{\small
\noindent
\(
1+\frac{A[1,0,0,0]}{q}+\frac{A[1,1,0,0]}{q}+\frac{A[1,1,2,0]}{q}
+\frac{A[1,1,2,2]}{q}+\frac{A[1,2,2,0]}{q}+\frac{A[1,2,2,2]}{q}
+\frac{A[1,2,4,2]}{q}+\frac{A[1,3,4,2]}{q}+\frac{A[2,2,2,0]}{q^2}
+\frac{A[2,2,2,2]}{q^2}+\frac{A[2,2,4,2]}{q^2}+\frac{A[2,4,4,2]}{q^2}
+\frac{A[2,4,6,2]}{q^2}+\frac{A[2,4,6,4]}{q^2}+\frac{A[3,3,4,2]}{q^2}
+\frac{A[3,4,4,2]}{q^2}+\frac{A[3,4,6,2]}{q^2}+\frac{A[3,4,6,4]}{q^2}
+\frac{A[3,5,6,2]}{q^2}+\frac{A[3,5,6,4]}{q^2}+\frac{A[3,5,8,4]}{q^2}
+\frac{A[3,6,8,4]}{q^2}+\frac{A[4,6,8,4]}{q^2}.
\)
}
\renewcommand{\baselinestretch}{1.2} 
\smallskip

$\tilde{\E}^\dag_2=$
\renewcommand{\baselinestretch}{0.5} 
{\small
\noindent
\(
1+\frac{A[0,1,0,0]}{q}+\frac{A[0,1,2,0]}{q}+\frac{A[0,1,2,2]}{q}
+\frac{A[0,2,2,0]}{q^2}+\frac{A[0,2,2,2]}{q^2}+\frac{A[0,2,4,2]}{q^2}
+\frac{A[0,3,4,2]}{q^3}+\frac{A[1,1,0,0]}{q}+\frac{A[1,1,2,0]}{q}
+\frac{A[1,1,2,2]}{q}+\frac{A[1,3,2,0]}{q^2}+\frac{A[1,3,2,2]}{q^2}
+\frac{A[1,3,4,0]}{q^2}+\frac{A[1,3,4,4]}{q^2}+\frac{A[1,3,6,2]}{q^2}
+\frac{A[1,3,6,4]}{q^2}+\frac{A[1,5,6,2]}{q^3}+\frac{A[1,5,6,4]}{q^3}
+\frac{A[1,5,8,4]}{q^3}+\frac{A[2,2,2,0]}{q^2}+\frac{A[2,2,2,2]}{q^2}
+\frac{A[2,2,4,2]}{q^2}+\frac{A[2,3,2,0]}{q^2}+\frac{A[2,3,2,2]}{q^2}
+\frac{A[2,3,4,0]}{q^2}+\frac{A[2,3,4,4]}{q^2}+\frac{A[2,3,6,2]}{q^2}
+\frac{A[2,3,6,4]}{q^2}+\frac{A[2,4,4,0]}{q^2}+\frac{A[2,4,4,4]}{q^2}
+\frac{A[2,4,8,4]}{q^2}+\frac{A[2,5,4,2]}{q^3}+\frac{A[2,5,8,2]}{q^3}
+\frac{A[2,5,8,6]}{q^3}+\frac{A[2,6,6,2]}{q^3}+\frac{A[2,6,6,4]}{q^3}
+\frac{A[2,6,8,2]}{q^3}+\frac{A[2,6,8,6]}{q^3}+\frac{A[2,6,10,4]}{q^3}
+\frac{A[2,6,10,6]}{q^3}+\frac{A[2,7,8,4]}{q^3}
+\frac{A[2,7,10,4]}{q^3}
+\frac{A[2,7,10,6]}{q^3}+\frac{A[3,3,4,2]}{q^3}+\frac{A[3,5,4,2]}{q^3}
+\frac{A[3,5,8,2]}{q^3}+\frac{A[3,5,8,6]}{q^3}+\frac{A[3,7,8,2]}{q^3}
+\frac{A[3,7,8,6]}{q^3}+\frac{A[3,7,12,6]}{q^3}
+\frac{A[3,9,12,6]}{q^3}
+\frac{A[4,5,6,2]}{q^3}+\frac{A[4,5,6,4]}{q^3}+\frac{A[4,5,8,4]}{q^3}
+\frac{A[4,6,6,2]}{q^3}+\frac{A[4,6,6,4]}{q^3}+\frac{A[4,6,8,2]}{q^3}
+\frac{A[4,6,8,6]}{q^3}+\frac{A[4,6,10,4]}{q^3}
+\frac{A[4,6,10,6]}{q^3}
+\frac{A[4,7,8,2]}{q^3}+\frac{A[4,7,8,6]}{q^3}+\frac{A[4,7,12,6]}{q^3}
+\frac{A[4,8,8,4]}{q^4}+\frac{A[4,8,12,4]}{q^4}
+\frac{A[4,8,12,8]}{q^4}
+\frac{A[4,9,10,4]}{q^4}+\frac{A[4,9,10,6]}{q^4}
+\frac{A[4,9,12,4]}{q^4}
+\frac{A[4,9,12,8]}{q^4}+\frac{A[4,9,14,6]}{q^4}
+\frac{A[4,9,14,8]}{q^4}
+\frac{A[4,10,12,6]}{q^4}+\frac{A[4,10,14,6]}{q^4}
+\frac{A[4,10,14,8]}{q^4}+\frac{A[5,7,8,4]}{q^3}
+\frac{A[5,7,10,4]}{q^3}
+\frac{A[5,7,10,6]}{q^3}+\frac{A[5,9,10,4]}{q^4}
+\frac{A[5,9,10,6]}{q^4}
+\frac{A[5,9,12,4]}{q^4}+\frac{A[5,9,12,8]}{q^4}
+\frac{A[5,9,14,6]}{q^4}
+\frac{A[5,9,14,8]}{q^4}+\frac{A[5,11,14,6]}{q^5}
+\frac{A[5,11,14,8]}{q^5}+\frac{A[5,11,16,8]}{q^5}
+\frac{A[6,9,12,6]}{q^3}+\frac{A[6,10,12,6]}{q^4}
+\frac{A[6,10,14,6]}{q^4}+\frac{A[6,10,14,8]}{q^4}
+\frac{A[6,11,14,6]}{q^5}+\frac{A[6,11,14,8]}{q^5}
+\frac{A[6,11,16,8]}{q^5}+\frac{A[6,12,16,8]}{q^6}.
\)
}
\renewcommand{\baselinestretch}{1.2} 
\smallskip

$\tilde{\E}^\dag_3=$
\renewcommand{\baselinestretch}{0.5} 
{\small
\noindent
\(
1+\frac{A[0,0,1,0]}{q}+\frac{A[0,0,1,1]}{q}+\frac{A[0,1,1,0]}{q}
+\frac{A[0,1,1,1]}{q}+\frac{A[0,1,2,0]}{q}+\frac{A[0,1,2,2]}{q}
+\frac{A[0,1,3,1]}{q^2}+\frac{A[0,1,3,2]}{q^2}+\frac{A[0,2,3,1]}{q^2}
+\frac{A[0,2,3,2]}{q^2}+\frac{A[0,2,4,2]}{q^2}+\frac{A[1,1,1,0]}{q}
+\frac{A[1,1,1,1]}{q}+\frac{A[1,1,2,0]}{q}+\frac{A[1,1,2,2]}{q}
+\frac{A[1,1,3,1]}{q^2}+\frac{A[1,1,3,2]}{q^2}+\frac{A[1,2,2,0]}{q}
+\frac{A[1,2,2,2]}{q}+\frac{A[1,2,3,0]}{q^2}+\frac{A[1,2,3,3]}{q^2}
+\frac{A[1,2,5,2]}{q^2}+\frac{A[1,2,5,3]}{q^2}+\frac{A[1,3,3,1]}{q^2}
+\frac{A[1,3,3,2]}{q^2}+\frac{A[1,3,5,1]}{q^2}+\frac{A[1,3,5,4]}{q^2}
+\frac{A[1,3,6,2]}{q^2}+\frac{A[1,3,6,4]}{q^2}+\frac{A[1,4,5,2]}{q^2}
+\frac{A[1,4,5,3]}{q^2}+\frac{A[1,4,6,2]}{q^2}+\frac{A[1,4,6,4]}{q^2}
+\frac{A[1,4,7,3]}{q^3}+\frac{A[1,4,7,4]}{q^3}+\frac{A[2,2,3,1]}{q^2}
+\frac{A[2,2,3,2]}{q^2}+\frac{A[2,2,4,2]}{q^2}+\frac{A[2,3,3,1]}{q^2}
+\frac{A[2,3,3,2]}{q^2}+\frac{A[2,3,5,1]}{q^2}+\frac{A[2,3,5,4]}{q^2}
+\frac{A[2,3,6,2]}{q^2}+\frac{A[2,3,6,4]}{q^2}+\frac{A[2,4,4,2]}{q^2}
+\frac{A[2,4,5,1]}{q^2}+\frac{A[2,4,5,4]}{q^2}+\frac{A[2,4,7,2]}{q^3}
+\frac{A[2,4,7,5]}{q^3}+\frac{A[2,4,8,4]}{q^2}+\frac{A[2,5,6,2]}{q^2}
+\frac{A[2,5,6,4]}{q^2}+\frac{A[2,5,7,2]}{q^3}+\frac{A[2,5,7,5]}{q^3}
+\frac{A[2,5,9,4]}{q^3}+\frac{A[2,5,9,5]}{q^3}+\frac{A[2,6,8,4]}{q^2}
+\frac{A[2,6,9,4]}{q^3}+\frac{A[2,6,9,5]}{q^3}+\frac{A[3,4,5,2]}{q^2}
+\frac{A[3,4,5,3]}{q^2}+\frac{A[3,4,6,2]}{q^2}+\frac{A[3,4,6,4]}{q^2}
+\frac{A[3,4,7,3]}{q^3}+\frac{A[3,4,7,4]}{q^3}+\frac{A[3,5,6,2]}{q^2}
+\frac{A[3,5,6,4]}{q^2}+\frac{A[3,5,7,2]}{q^3}+\frac{A[3,5,7,5]}{q^3}
+\frac{A[3,5,9,4]}{q^3}+\frac{A[3,5,9,5]}{q^3}+\frac{A[3,6,7,3]}{q^3}
+\frac{A[3,6,7,4]}{q^3}+\frac{A[3,6,9,3]}{q^3}+\frac{A[3,6,9,6]}{q^3}
+\frac{A[3,6,10,4]}{q^3}+\frac{A[3,6,10,6]}{q^3}
+\frac{A[3,7,9,4]}{q^3}
+\frac{A[3,7,9,5]}{q^3}+\frac{A[3,7,10,4]}{q^3}
+\frac{A[3,7,10,6]}{q^3}
+\frac{A[3,7,11,5]}{q^4}+\frac{A[3,7,11,6]}{q^4}
+\frac{A[4,6,8,4]}{q^2}
+\frac{A[4,6,9,4]}{q^3}+\frac{A[4,6,9,5]}{q^3}+\frac{A[4,7,9,4]}{q^3}
+\frac{A[4,7,9,5]}{q^3}+\frac{A[4,7,10,4]}{q^3}
+\frac{A[4,7,10,6]}{q^3}
+\frac{A[4,7,11,5]}{q^4}+\frac{A[4,7,11,6]}{q^4}
+\frac{A[4,8,11,5]}{q^4}
+\frac{A[4,8,11,6]}{q^4}+\frac{A[4,8,12,6]}{q^4}.
\)
}
\renewcommand{\baselinestretch}{1.2} 
\smallskip

$\tilde{\E}^\dag_4=$
\renewcommand{\baselinestretch}{0.5} 
{\small
\noindent
\(
1+\frac{A[0,0,0,1]}{q}+\frac{A[0,0,1,1]}{q}+\frac{A[0,1,1,1]}{q}
+\frac{A[0,1,2,1]}{q}+\frac{A[0,1,2,2]}{q}+\frac{A[1,1,1,1]}{q}
+\frac{A[1,1,2,1]}{q}+\frac{A[1,1,2,2]}{q}+\frac{A[1,2,2,1]}{q}
+\frac{A[1,2,2,2]}{q}+\frac{A[1,2,3,1]}{q}+\frac{A[1,2,3,3]}{q^2}
+\frac{A[1,2,4,2]}{q}+\frac{A[1,2,4,3]}{q^2}+\frac{A[1,3,4,2]}{q}
+\frac{A[1,3,4,3]}{q^2}+\frac{A[1,3,5,3]}{q^2}+\frac{A[2,3,4,2]}{q}
+\frac{A[2,3,4,3]}{q^2}+\frac{A[2,3,5,3]}{q^2}+\frac{A[2,4,5,3]}{q^2}
+\frac{A[2,4,6,3]}{q^2}+\frac{A[2,4,6,4]}{q^2}.
\)
}
\renewcommand{\baselinestretch}{1.2} 
\medskip

{\bf Acknowledgments.}
The first author thanks IHES for the invitation. 
We thank Daniel Orr for useful discussions
on the positivity conjectures from \cite{CO1} and Dmitry Timashev 
for explanations on the Vinberg filtration.

The work of Ivan Cherednik was partially supported by 
{\small\em NSF grant DMS--1101535}. 

The work of Evgeny Feigin was partially supported
by {\small\em the Russian President Grant MK-3312.2012.1, 
{\rm by} the Dynasty Foundation,
{\rm by} the AG Laboratory HSE, RF government grant, 
ag. 11.G34.31.0023, 
{\rm by} the RFBR grants 12-01-00070, 12-01-00944, 12-01-33101  
{\rm and by} the Russian Ministry of Education and Science under the
grant 2012-1.1-12-000-1011-016.}
This study comprises research fundings from the 
{\small\em 
"Representation Theory in Geometry and in Mathematical Physics" 
carried out within the
National Research University Higher School of Economics' 
Academic Fund Program in 2012, grant No 12-05-0014}.
This study was carried out within the {\small\em National Research 
University Higher School of Economics
Academic Fund Program in 2012-2013, research grant No. 11-01-0017}.

\comment{
{\em Sketch of the correction to Conjecture (only in $V_\om$)}.
For $D_4$. 
Let $b=\{0, 0, 1, 1\}$ in terms of the fundamental weights.
Then $b-w_0(b)=11+12+10$;
$(\om_3-w_0(\om_3))+(\om_4-w_0(\om_4))= (12+3)+(12+4).$
Thus, to save the conjecture,
$f_{11}f_{12}f_{10}(vac)$ must be $0$ in $V_b$ and
$f_{12}^2 f_3 f_4 (vac)\neq 0$ there.

We consider $V_1\otimes V_2$ for
$V_1=V_{\om_3}, V_2=V_{\om_4}$ with the highest vectors
$v_{1,2}$.

We expect that 
$f_{12}f_3(vac_1)\otimes f_{12}f_4(vac_2)\neq 0$
in $V_b$. We then need to check that,
for instance, $u=f_{12}f_{4}\otimes f_{12}f_4(v_1\otimes v_2)$
is zero in $V_b$. We can use that the map $V_1\otimes V_2\to V_b$
becomes an isomorphism from the linear span of 
$\C w(v_1\otimes v_2)$
and $\C w(v_b)$ for all $w\in W$ (the multiplicity one theorem
in $V_b$). The kernel of this homomorphism is the greatest
such $\mathfrak{g}$-submodule. Why $\mathfrak{g}(u)$ is a proper 
submodule?

Similarly, why the image of $f_{11}f_{12}f_{10}(v_1\otimes v_2)$
is zero? We need to check, for instance, that
the image of $f_{11}\otimes f_{12}f_{10}(v_1\otimes v_2)$
is zero.  
}

\bibliographystyle{unsrt}

\end{document}